\documentclass[journal]{IEEEtai}

\usepackage[colorlinks,urlcolor=blue,linkcolor=blue,citecolor=blue]{hyperref}
\usepackage{color,array}
\usepackage{graphicx}
\usepackage[super]{nth}
\usepackage{amsbsy, amsmath, amsfonts}
\usepackage{amssymb}
\usepackage{mathtools}
\usepackage{stmaryrd}
\usepackage{commath}
\usepackage{physics}
\usepackage{esvect}
\usepackage{upgreek}
\usepackage[hang,flushmargin]{footmisc}
\usepackage{caption}
\usepackage{subcaption}
\usepackage{xcolor}
\usepackage{bbding}

\usepackage{amsthm}
\theoremstyle{definition}
\newtheorem{definition}{Definition}

\usepackage{tikz}
\newcommand*\circled[1]{\tikz[baseline=(char.base)]{
            \node[shape=circle,draw,inner sep=1pt] (char) {#1};}}

\usepackage{pifont}
\usepackage{tabularx}
\usepackage{makecell}
\usepackage{makecell}
\usepackage{csquotes}
\usepackage{url}
\usepackage{multirow}

\usepackage{algpseudocode}
\usepackage{algorithm}
\usepackage{algpseudocode}

\newcommand\inv[1]{#1\raisebox{1.15ex}{$\scriptscriptstyle-\!1$}}

\DeclareMathOperator{\E}{\mathbb{E}}
\DeclareMathOperator{\V}{\mathbb{V}}

\DeclareMathOperator*{\argmin}{arg\,min}

\makeatletter
\newcommand{\xMapsto}[2][]{\ext@arrow 0599{\Mapstofill@}{#1}{#2}}
\def\Mapstofill@{\arrowfill@{\Mapstochar\Relbar}\Relbar\Rightarrow}
\makeatother

\begin{document}
\title{Noise-aware Physics-informed Machine Learning for Robust PDE Discovery}

\author{Pongpisit Thanasutives, Takashi Morita, Masayuki Numao, and Ken-ichi Fukui
\thanks{Pongpisit Thanasutives is with Graduate School of Information Science and Technology, Osaka University, Japan (e-mail: thanasutives@ai.sanken.osaka-u.ac.jp).}
\thanks{Takashi Morita, Masayuki Numao, and Ken-ichi Fukui are with Osaka University, Japan (e-mail: \{t-morita, numao, fukui\}@ai.sanken.osaka-u.ac.jp).}
}

\maketitle

\begin{abstract}
This work is concerned with discovering the governing partial differential equation (PDE) of a physical system. Existing methods have demonstrated the PDE identification from finite observations but failed to maintain satisfying results against noisy data, partly owing to suboptimal estimated derivatives and found PDE coefficients. We address the issues by introducing a noise-aware physics-informed machine learning (nPIML) framework to discover the governing PDE from data following arbitrary distributions. We propose training a couple of neural networks, namely solver and preselector, in a multi-task learning paradigm, which yields important scores of basis candidates that constitute the hidden physical constraint. After they are jointly trained, the solver network estimates potential candidates, e.g., partial derivatives, for the sparse regression algorithm to initially unveil the most likely parsimonious PDE, decided according to the information criterion. We also propose the denoising physics-informed neural networks (dPINNs), based on Discrete Fourier Transform (DFT), to deliver a set of the optimal finetuned PDE coefficients respecting the noise-reduced variables. The denoising PINNs are structured into forefront projection networks and a PINN, by which the formerly learned solver initializes. Our extensive experiments on five canonical PDEs affirm that the proposed framework presents a robust and interpretable approach for PDE discovery, applicable to a wide range of systems, possibly complicated by noise.
\end{abstract}


\section{Introduction} \label{introduction}
Data-driven discovery has recently gained popularity due to its flexibility and satisfactory accuracy in uncovering the hidden underlying partial differential equation (PDE) of a dynamical system with less required domain knowledge. Applying sparse regression-based approaches to a library of the target variable and its partial derivative candidates is a promising method for discovering a parsimonious model purely out of observational data. A few of such previous attempts were, for instance, sequential threshold ridge regression (STRidge) \cite{rudy2017data}, $L_{1}$-regularized sparse optimization  \cite{schaeffer2017learning} based on the least absolute shrinkage and selection operator (LASSO) \cite{tibshirani1996regression}, and sparse Bayesian regression \cite{zhang2018robust}.

Since partial derivatives are treated as the vital input features, inaccurate estimation of the derivatives, primarily the high-order ones, using numerical differentiation, such as finite difference, whose performance drops when facing sparse corrupted data, can poorly affect the discovered results. This paper utilizes automatic differentiation (AD) \cite{baydin2017automatic} on a neural network that we refer to as the solver to be an alternative approach, as formerly suggested by \cite{raissi2019physics}. AD allows derivative computation given a mere implementation of the hypothesis function; therefore, the method does not suffer from truncation error, mitigating the numerical imprecision of computing high-order derivatives.

Although the utilization of neural networks is not restricted by the assumption of particular input distributions, the solver that learns just by correcting prediction errors may be prone to overfitting to the finite observations and inadequate for capturing the proper PDE solution, especially when encountering sparse measurements. This troublesome motivates us to formulate the solver network with weak physics-informed regularization, maintaining the prediction performance while respecting an implicit form of the governing physical law. Specifically new in this work, we propose multi-task training with a preselector neural network that promotes sparsity based on an interpretable self-gating mechanism to alleviate the issue. The preselector learns the system's estimated evolution (produced by the solver) from the spatial derivatives and other features to represent the hidden parsimonious PDE. Furthermore, we present a workable way of using the trained preselector's feature importance to encourage selecting the expressive candidates that derive a non-overfitting PDE.

Once both the networks are converged using a multi-task learning procedure and the library of potential (nonlinear) terms is prepared, we then apply a form of sparse linear regression algorithms, e.g., STRidge \cite{rudy2017data}, to the discretized domain of interest. Nonetheless, the proper selection of the regularization hyperparameters regarding the sparse linear model can be problematic. While the true underlying PDE remains unknown, solely cross-validating equations together with the Pareto analysis based on one fixed-valued regularization hyperparameter may still yield insignificant, probably wrong results, especially in a small data regime. Thus, as an additional consideration, the initial discovered PDE is encouraged to include the candidates whose importance is greater than a threshold defined within the proposed preselector network's $L_{0}$-regularized self-gating mechanism. After the cooperative learning, among the expected PDEs formable using the threshold-passing basis candidates, the parsimonious but informative PDEs are preferred, i.e., having sufficiently low Bayesian information criterion (BIC) \cite{schwarz1978estimating} or Akaike information criterion (AIC) \cite{akaike1998information}.

At this point, the sparse learning algorithm has yielded a guess of the hidden governing PDE, referred to as the initial discovered PDE. However, propagating error is woefully inevitable since the sparse regression is separated from the candidate library preparation step. This consequently causes the initial PDE to not be at its optimum concerning the given input data. To achieve the most-favorable PDE, we parameterize all the discovered coefficients as the gradient-based learnable parameters of the physics-informed solver network that is finetuned such that its output approximates the target variable while concurrently respecting the most-relevant underlying PDE as per the core proposal of physics-informed learning \cite{raissi2019physics, karniadakis2021physics}. Remark that, without an appropriate initialization of targeting PDE coefficients, training a physics-informed neural network (PINN) \cite{raissi2019physics} may be a task that could be developed further by, for example, multi-task learning \cite{thanasutives2021adversarial} or sinusoidal feature mapping \cite{wong2021learning}, even though the actual governing function is presumably known beforehand.

In a practical scenario where noise may disturb both the independent and dependent variables, the optimization process of PINN is perturbed; thus, attaining a local optimum set of coefficients is spaced out from the ground truth. The previous work, abbreviated as DLrSR \cite{li2020robust}, tackled the difficulty via low-rank matrix factorization solved by robust PCA \cite{candes2011robust}, neglecting the assumably sparse noise and utilizing the low-rank data. Nevertheless, if the sparse noise presumption does not hold, the method can be impotent for various situations. To mitigate the issue, we introduce denoising layers based on precomputed Discrete Fourier Transform (DFT) to the vanilla PINN, optimizing the solver-founded PDE. The denoising layers filter out the frequency components of the input signal, whose power is less than a predefined threshold, then obtain contaminated noise by taking the difference between the original and reconstructed signal. The extracted noises are projected using projection neural networks to perturb backwardly or denoise the noisy measurements with the appropriate intensities, then reconstruct the noise-reduced dataset. This paper coins a PINN attached to the proposed denoising mechanism as denoising PINNs (dPINNs). Ultimately, after the dPINNs' learning, the converged parameters regarding all effective coefficients are treated as the end results.

Experimental results from 5 canonical models, including 3 ordinary PDEs and 2 complex-valued PDEs, reveal that the proposed framework outperforms the state-of-the-art sparse regression methods in noiseless and noisy datasets. As a proof of concept distinct from prior works, we conduct investigations on learning from noisy independent variables, e.g., polluted spatial and temporal variables, which are relevant to GPS coordinate measurements \cite{ranacher2016gps} and manual timing in physical experiments \cite{faux2019manual}.

We summarize our main contributions as follows:
\begin{itemize}
  \item We introduce the multi-task learning with the preselector network to impose the weak physical constraint which is calculable without labeled supervision.
  \item We introduce an utilization of the preselector's perceived feature importance scores to bring an auxiliary view to the candidate selection, addressing the fundamental sensitivity problem of finding the right sparsity-promoting regularization on the sparse regression-based method.
  \item We introduce denoising physics-informed neural networks (dPINNs) based on DFT and the projection networks to handle both noisy independent and dependent variables.
\end{itemize}
\section{Method: Noise-aware Physics-informed Machine Learning (nPIML) Framework} \label{method}
\subsection{Problem Formulation and Overview} \label{problem_formulation}
We consider the following general form of nonlinear PDE in the dynamical system perspective:
\begin{equation} \label{eq:dynamical_system}
u_{t}=\mathcal{N}_{\xi}[\Theta];\quad \Theta=\mqty[u & u_{x} & u_{xx} & \cdots & x].
\end{equation}
\noindent \(\mathcal{N}_{\xi}\) is the governing function parameterized by the vector of coefficients \(\xi\). The function depends on \(\Theta\), which may consist of the spatial variable \(x\), the derivatives and any indispensable features. In regards to   \(\mathcal{N}_{\xi}\), \(\Theta\) is the smallest possible, merely composed of the necessary terms. \(u\) is the dependent PDE solution, observed with the space-time matrix \((x, t)\).

\begin{figure*}
	\centering
	\includegraphics[trim={0 1.25cm 0 0}, clip, width=0.9\linewidth, height=7cm]{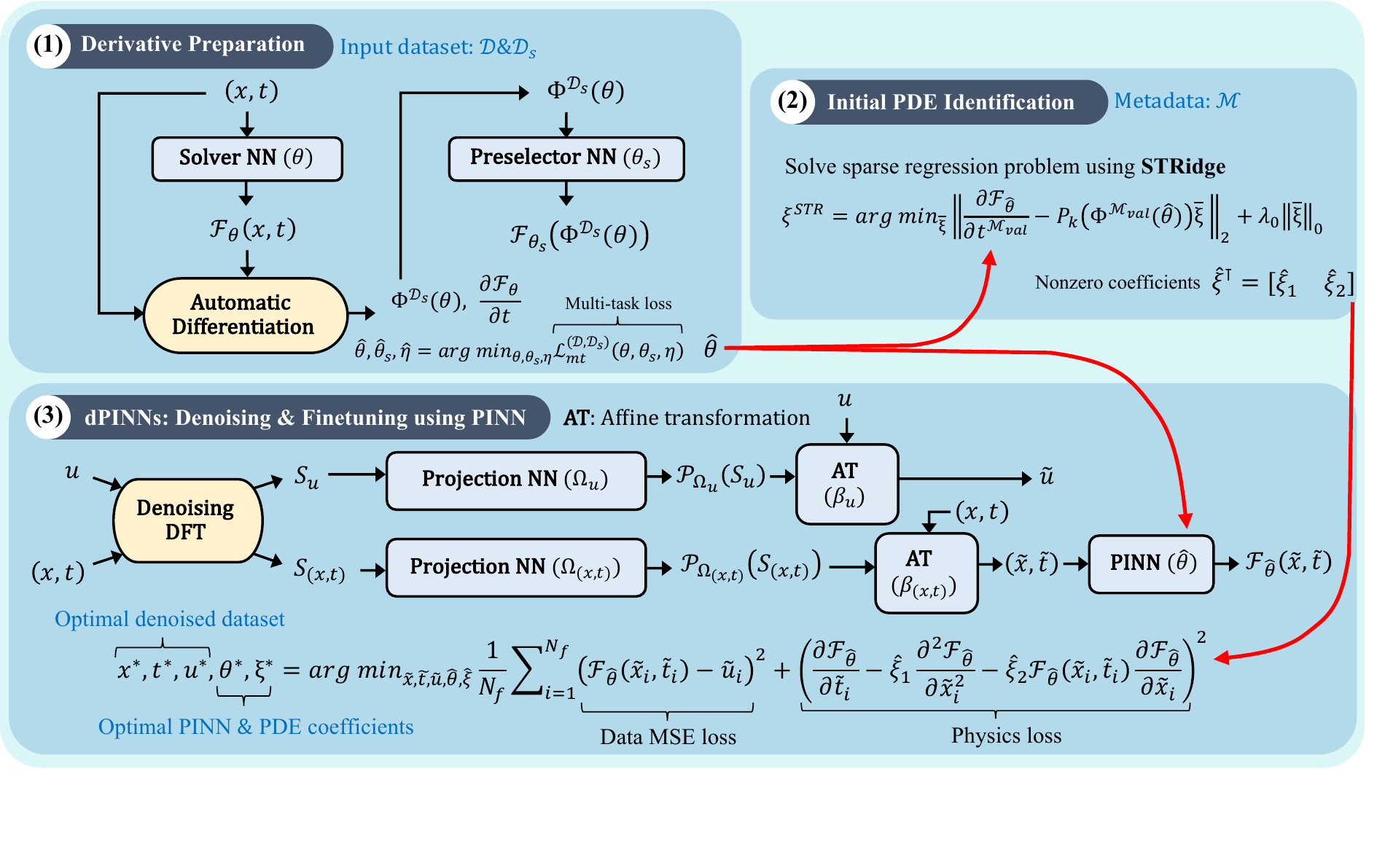}
	\caption{Exemplary discovery scheme of the proposed \textbf{noise-aware Physics-informed Machine Learning (nPIML) framework}: \textbf{(1)} Physics-regularized derivative preparation by multi-task learning of the solver and preselector. \textbf{(2)} Initial identification of the hidden PDE by STRidge. \textbf{(3)} Applying the denoising DFT to \((x, t)\)\&\(u\) then finetuning the initial PDE coefficients on the denoised variables, using PINN.}
	\label{fig:overview_framework}
\end{figure*}

Fig. \ref{fig:overview_framework}. conceptualizes the three principal procedures for uncovering \(\xi\) preferably in a low-dimensional space by walking through an exemplar of discovering Burgers' PDE \cite{basdevant1986spectral}. Step \circled{\textbf{(1)}}, we numerically equivalizes \(u\) and \(\mathcal{N}_{\xi}[\Theta]\) to the solver and preselector neural network outputs \(\mathcal{F}_{\theta}(x, t)\) and \(\mathcal{F}_{\theta_{s}}(\Phi^{\mathcal{D}_{s}}(\theta))\). \(\Phi^{\mathcal{D}_{s}}(\theta) \in \mathbb{C}^{(N_f+N_r) \times C}\) is the library of \(C\) linearly independent atomic/basis candidates from which the preselector learns to embed physics by inferring the system evolution. The candidates are evaluated on a set \(\mathcal{D}_{s}=\set{(x_{i}, t_{i})^{N_f+N_r}_{i=1}}\). Step \circled{\textbf{(2)}}, the well-fitted networks, \(\hat{\theta}\) and \(\hat{\theta}_{s}\), put together a larger library of potential \(k\)-degree polynomial features \(P_{k}(\Phi^{\mathcal{M}_{val}}(\hat{\theta}))\) of which an initial analytical expression of Burgers' PDE, worked out approximately by STRidge \cite{rudy2017data}, is made. Step \circled{\textbf{(3)}}, \(\hat{\theta}\) is henceforth transferred to the PINN that is optimally finetuned with the PDE, initialized by nonzero coefficients \(\hat{\xi}\), on the denoised variables \(\Tilde{x}\), \(\Tilde{t}\) and \(\Tilde{u}\), offered by the projection networks \(\mathcal{P}_{\Omega_{(x, t)}}\) and \(\mathcal{P}_{\Omega_{u}}\). The noise-reduction mechanism functions as a series of affine transformations, controlled by \(\beta_{(x, t)}\) and \(\beta_{u}\), of the dataset with the projected noises \(\mathcal{P}_{\Omega_{(x, t)}}(S_{(x, t)})\) and \(\mathcal{P}_{\Omega_{u}}(S_u)\), after applying frequency-based denoising DFT. The mathematical derivation of the relevant variables are elaborated more in \ref{derivative_preparation}, \ref{initial_pde_identification} and \ref{dPINNs}.

\subsection{Derivative Preparation} \label{derivative_preparation}
Concerning \circled{\textbf{(1)}} of Fig. \ref{fig:overview_framework}, we utilize the solver (\(\mathcal{F}_{\theta}\)) and preselector (\(\mathcal{F}_{\theta_{s}}\)) networks, which are jointly trained for the solver network to be weakly physics-constrained. Facilitating the co-training, the solver network is pretrained on the dataset \(\mathcal{D}=\set{(x_{i}, t_{i}, u_{i})^{N_f}_{i=1}}\) to approximate the mapping function. Therefore, the partial derivative candidate values are assured of becoming close to the valid values. At the pretraining stage, the solver network minimizes the mean square error (MSE)
\begin{equation} \label{eq:loss_sup}
\mathcal{L}^{\mathcal{D}}_{sup}(\theta) = \frac{1}{N_{f}}\sum^{N_{f}}_{i=1}(\mathcal{F}_{\theta}(x_{i}, t_{i})-u_{i})^{2};\, (x_{i}, t_{i}, u_{i}) \in \mathcal{D},
\end{equation}
\noindent where \(N_f\) is the number of labeled subsamples. If \(u\) is complex-valued, the sum of the MSEs from the real and imaginary parts is taken as the supervised loss function. Since the futile search over infinitely feasible \(\Phi^{\mathcal{D}_{s}}(\theta)\) setups would be intractable, we instead build an overcomplete candidate library given to the preselector network for deciding the informative set of features by minimizing
\begin{equation} \label{eq:loss_unsup}
\begin{aligned}
\mathcal{L}^{\mathcal{D}_{s}}_{unsup}(\theta, \theta_{s}) &= \frac{1}{N_{f}+N_{r}}\sum^{N_{f}+N_{r}}_{i=1}(\frac{\partial\mathcal{F}_{\theta}}{\partial t_{i}}-\mathcal{F}_{\theta_{s}}(\Phi^{\mathcal{D}_{s}}_{i}(\theta)))^{2};\\
\Phi^{\mathcal{D}_{s}}_{i}(\theta) &= \mqty[\mathcal{F}_{\theta}(x_{i}, t_{i}) & \frac{\partial\mathcal{F}_{\theta}}{\partial x_{i}} & \frac{\partial^{2}\mathcal{F}_{\theta}}{\partial x^{2}_{i}} & \cdots &x_{i}],
\end{aligned}
\end{equation}
\noindent where \(N_r\) is the number of unsupervised subsamples within the domain that disjoints the supervised set \(\mathcal{D}\). We attain \(\mathcal{D}_{s}\) by fusing up the spatio-temporal measurements without supervision. Each derivative term's input is usually omitted for notational convenience. Inspired by the assumption that low-order partial derivatives are commonly included more than the higher ones, we embed the thresholded self-gated mechanism, parameterized by \(W^b\), to the preselector forward pass, emphasizing the priority of simple models as follows:
\begin{equation}
\begin{aligned}
\mathcal{F}_{\theta_s}(\Phi^{\mathcal{D}_{s}}(\theta)) &= \mathcal{F}_{\theta^{r}_s}(\mathcal{F}_{W^b}(\Phi^{\mathcal{D}_{s}}(\theta))),\\
\mathcal{F}_{W^b}(\Phi^{\mathcal{D}_{s}}(\theta)) &= \Phi^{\mathcal{D}_{s}}(\theta) \odot \mathcal{A}^{\mathcal{T}}(\Phi^{\mathcal{D}_{s}}(\theta), W^b),\\
\mathcal{A}^{\mathcal{T}}_{j}(\Phi^{\mathcal{D}_{s}}(\theta), W^b) &= \max(\mathcal{A}_j(\Phi^{\mathcal{D}_{s}}(\theta), W^b)-\mathcal{T}, 0),\\
\mathcal{A}_{j}(\Phi^{\mathcal{D}_{s}}(\theta), W^b) &= \frac{\sum^{N_f+N_r}_{i=1}\sigma(\sum^{C}_{k=1}\Phi^{\mathcal{D}_{s}}_{ik}(\theta)W_{kj}+b_{j})}{N_f+N_r}.
\end{aligned}
\end{equation}
\noindent \(\odot\) refers to Hadamard product (broadcast multiplication). \(\mathcal{A}^{\mathcal{T}}(\Phi^{\mathcal{D}_{s}}(\theta), W^b)\) is interpreted as the thresholded vector-valued feature importance the preselector perceive. The self-gated mechanism utilizes the activation function \(\sigma\) to compute the expected importance of each candidate in terms of (unnormalized) probability across \(N_f+N_r\) samples. Note that we only consider the real part of \(\Phi^{\mathcal{D}_{s}}(\theta)W + b\) in the case of complex-valued PDEs. \(\mathcal{T}\) is a threshold for allowing the effective basis candidates. The threshold is initialized to be surely less than the minimal candidate importance, specifically we set \(\mathcal{T} = \kappa \min_{j}\mathcal{A}^{(1)}_{j}(\Phi^{\mathcal{D}_{s}}(\theta), W^b)\), where \(0<\kappa<1\), before the first joint gradient update, denoted by the superscript \((1)\). The parameter \(W^b\) consists of \(W \in \mathbb{C}^{C \times C}\) and \(b \in \mathbb{C}^{1 \times C}\) (weights and biases of the linear layer), serving as the share of the preselector's parameters:
\begin{equation}
\theta_s = (W^b, \theta^{r}_s);\, \mathcal{F}_{\theta_s} =  \mathcal{F}_{\theta^{r}_s} \circ \mathcal{F}_{W^b}.
\end{equation}
\noindent Excluding \(W^b\), the rest of the preselector network's parameters get referred to as \(\theta^{r}_s\). We devise \(R^{\mathcal{D}_{s}}(\theta, W^b)\) as a $L_0$-regularization on \(\mathcal{A}^{\mathcal{T}}\) for selecting the expressive subset with priority to lower-order candidates in favor of Occam's razor principle. The regularization, encouraging the sparse and simple preselector learned representations, reads
\begin{equation} \label{eq:reg}
\begin{aligned}
    R^{\mathcal{D}_{s}}(\theta, W^b) &= \lambda_{1}(\norm{\mathcal{A}^{\mathcal{T}}(\Phi^{\mathcal{D}_{s}}(\theta), W^b)}_{0}\\
    &+\lambda_{2}\sum^{C}_{j=1}w_{j}\mathcal{A}^{\mathcal{T}}_{j}(\Phi^{\mathcal{D}_{s}}(\theta), W^b)).\\
\end{aligned}
\end{equation}
\noindent \(w\) is the weighting by derivative orders, directly applied to the feature importance. For instance, suppose that \(j^{\text{th}}\) basis candidate associates to the second-order derivative \(u_{xx}\). Then we have \(w_{j}=2\). For nonderivative terms, we assign \(w_{j}=1\). \(\lambda_{1}\) is the parameter that controls the regularization intensity. \(\lambda_{2}\) closes the gap between the derivative orders such that the high-order derivatives are not always deselected. To practically minimize \(R^{\mathcal{D}_{s}}(\theta, W^{b})\) with \(\mathcal{L}^{\mathcal{D}_{s}}_{unsup}(\theta, \theta_{s})\) by a gradient-based optimizer, we have to overcome the obstacle that the \(L_{0}\) norm is not yet readily differentiable with respect to its input vector. Unlike how the gradient-free STRidge algorithm is executed, we require the smooth approximated function of \(L_{0}\) for achieving the thresholded feature importance. Adapted from SL0 algorithm \cite{mohimani2007fast}, we estimates
\begin{multline} \label{eq:l0_relaxation}
\norm{\mathcal{A}^{\mathcal{T}}(\Phi^{\mathcal{D}_{s}}(\theta), W^b)}_{0} \approx C\\ - \sum^{C}_{j=1}\exp(\frac{-(\mathcal{A}^{\mathcal{T}}_{j}(\Phi^{\mathcal{D}_{s}}(\theta), W^b))^{2}}{2(\eta\V(\mathcal{A}^{\mathcal{T}}(\Phi^{\mathcal{D}_{s}}(\theta), W^b)))^{2}}),
\end{multline}
\noindent where \(\V\) is the unbiased variance estimator over the \(C\) basis candidates. \(\eta\) determines the trade-off between the accuracy and smoothness: the smaller \(\eta\) gives the closer approximation, and the larger \(\eta\) gives the smoother approximation. \(\eta\) is initialized at \(1.0\) and learned with the gradients. We now denote the differentiable regularization function as \(R^{\mathcal{D}_{s}}_{\eta}(\theta, W^b)\). Combining (\ref{eq:loss_sup}), (\ref{eq:loss_unsup}), (\ref{eq:reg}) and (\ref{eq:l0_relaxation}), we view the multi-task learning of the weakly physics-informed solver and the coordinating simplicity-guided preselector inherently as the semi-supervised multi-objective optimization formulated as follows:
\begin{multline}\label{eq:loss_mtl}
\begin{aligned}
\hat{\theta}, \hat{\theta}_{s}, \hat{\eta} &= \argmin_{\theta, \theta_{s}, \eta} \mathcal{L}^{(\mathcal{D}, \mathcal{D}_{s})}_{mt}(\theta, \theta_{s}, \eta);\\
\mathcal{L}^{(\mathcal{D}, \mathcal{D}_{s})}_{mt}(\theta, \theta_{s}, \eta) &= MT(\mathcal{L}^{\mathcal{D}}_{sup}(\theta),
\end{aligned}
\\ \mathcal{L}^{\mathcal{D}_{s}}_{unsup}(\theta, \theta_{s})+R^{\mathcal{D}_{s}}_{\eta}(\theta, W^b)).
\end{multline}
\noindent The parameters of both networks are concurrently updated with the expectancy that the preselector network distills the hidden PDE function \(\mathcal{N}_{\xi}\), and informs physics back to the solver. \(MT\) is a function that reasonably manipulates learning by multiple losses, such as Uncert \cite{kendall2018multi} and PCGrad \cite{yu2020gradient}, which are shown to accelerate the PINN generalized performance \cite{thanasutives2021adversarial}. \textbf{Algorithm \ref{alg:1}} describes a relaxed approach that numerically minimizes the loss in (\ref{eq:loss_mtl}) until detected plateau; then, converging the solver network independently.

\begin{algorithm}
\caption{Multi-task learning for Derivative Preparation and Initial PDE Identification}
\label{alg:1}
\begin{algorithmic}[1]
\State \textbf{Goal:} To initially discover the governing function $\hat{\mathcal{N}}_{\hat{\xi}}$ based on the solver and preselector parameters $\hat{\theta}, \hat{\theta}_s$.
\State \textbf{Require:} Pretrained $\theta$ by (\ref{eq:loss_sup}) \& initialized $\theta_{s}$
\State Joint train\footnotemark[1] $\theta, \hat{\theta}_{s}, \hat{\eta} \leftarrow \argmin_{\theta, \theta_{s}, \eta} \mathcal{L}_{mt}^{(\mathcal{D}, \mathcal{D}_{s})}(\theta, \theta_{s}, \eta)$ \label{alg:1:joint_train}
\State Assign $I_{j} \leftarrow \mathcal{A}_j(\Phi^{\mathcal{D}_{s}}(\theta, \hat{W}^{\hat{b}}))-\mathcal{T}+\frac{1}{C}$ as the feature importance for each $j^{\text{th}}$ basis candidate
\State Converge the solver $\hat{\theta} \leftarrow \argmin_\theta\mathcal{L}^{\mathcal{D}}_{sup}(\theta)$
\State Build the candidate library on the metadata $\mathcal{M}$ by 
$\Phi^{\mathcal{M}}(\hat{\theta}) \leftarrow \mqty[\mathcal{F}_{\hat{\theta}}(x^{\mathcal{M}}, t^{\mathcal{M}}) &  \frac{\partial\mathcal{F}_{\hat{\theta}}}{\partial x^{\mathcal{M}}} &  \frac{\partial^{2}\mathcal{F}_{\hat{\theta}}}{\partial (x^{\mathcal{M}})^{2}} & \cdots & x^{\mathcal{M}}]$
\State Find $\hat{\mathcal{N}}_{\hat{\xi}}$ on $P_{k}(\Phi^{\mathcal{M}_{val}}(\hat{\theta}))$ using $\lambda_{STR}$-varied STRidge
\State \textbf{Return:} $\hat{\theta}, \hat{\theta}_s=(\hat{W}^{\hat{b}}, \hat{\theta}^{r}_s)$, $\hat{\eta}$ and $\hat{\mathcal{N}}_{\hat{\xi}}$
\end{algorithmic}
\hrulefill \\
\footnotesize \footnotemark[1]{After the joint training until empirical plateau, the learned preselector's parameters are regarded as $\hat{\theta}_{s}$. Converging the preselector could have been done, i.e., $\min_{\hat{\theta}^{r}_s} \mathcal{L}^{\mathcal{D}_{s}}_{unsup}(\theta, \hat{\theta}_{s})$, but did not to reduce the run time.}
\end{algorithm}

\subsection{Initial PDE identification} \label{initial_pde_identification}
Depicted by \circled{\textbf{(2)}} of Fig. \ref{fig:overview_framework}, we train STRidge \cite{rudy2017data} on top of the candidates and their polynomial features up to \(k\) degree: \(P_{k}(\Phi^{\mathcal{M}}(\hat{\theta}))\), which is evaluated on metadata \(\mathcal{M}\). For example, assume that \(k=2\), the unbiased interaction-only polynomial features of \(u,u_x\) and \(u_{xx}\) are formed as
\begin{equation}
P_{2}(\mqty[u & u_x & u_{xx}]) = \mqty[u & u_x & u_{xx} & uu_x & uu_{xx} & u_{x}u_{xx}].
\end{equation}
\noindent The metadata \(\mathcal{M}=\set{(x^{\mathcal{M}}_{i}, t^{\mathcal{M}}_{i})^{N_{\mathcal{M}}}_{i=1}}\) can be samples from a desired domain of interest, e.g., linearly discretized samples within a bounded rectangle domain are generated with the equal spaces as follows: \(\Delta x = \min_{i,j,(i \neq j)}\abs{x_i-x_j}\) and \(\Delta t = \min_{i,j,(i \neq j)}\abs{t_i-t_j}\). In fact, naively equating \(\forall i \in \set{1, 2, \dots, N_f}, (x^{\mathcal{M}}_{i}, t^{\mathcal{M}}_{i}) = (x_{i}, t_{i})\) is also viable for identifying the governing PDE as \(\hat{\mathcal{N}}_{\hat{\xi}}[P_{k}(\Phi^{\mathcal{M}_{val}}(\hat{\theta}))\mathcal{E}]\), where \(\hat{\xi}\) and \(\mathcal{E}\) are found by the following selection criterion:
\begin{equation} \label{eq:apply_stridge}
\begin{aligned}
\xi^{STR} &= \argmin_{\overline{\xi}} \norm{\frac{\partial \mathcal{F}_{\hat{\theta}}}{\partial t^{\mathcal{M}_{val}}} - P_{k}(\Phi^{\mathcal{M}_{val}}(\hat{\theta}))\overline{\xi}}_{2} + \lambda_{0}\norm{\overline{\xi}}_{0};\\
\lambda_{0} &= \mu\lambda_{STR}\varepsilon,\, E = \set{f_{i+1} \mid i \in \mathbb{N}_{\norm{\xi^{STR}}_{0}} \land \xi^{STR}_{f_{i+1}} \neq 0},\\
\hat{\xi} &= \mqty[\xi^{STR}_{f_1} & \cdots & \xi^{STR}_{f_{\abs{E}}}]^{\intercal},\, \mathcal{E} = \mqty[\boldsymbol{e}_{f_1} & \cdots & \boldsymbol{e}_{f_{\abs{E}}}].\\
\end{aligned}
\end{equation}
\noindent \(\varepsilon=\varepsilon(P_{k}(\Phi^{\mathcal{M}}(\hat{\theta})))\) is the significand of the conditional number (written in the scientific notation) of the candidate library. \(\mathcal{M}_{val}\) is a 20\% of the full \(\mathcal{M}\). For a tolerance \(tol\), \(\overline{\xi}\) is estimated by solving a relaxed \(\lambda_{STR}\)-regularized ridge regression problem on \(P_{k}(\Phi^{\mathcal{M}}(\hat{\theta}))\), whose polynomial candidate is normalized by its \(L_{2}\)-norm unless noted otherwise, with hard thresholding. To attain \(\xi^{STR}\), \(tol\) is iteratively refined with respect to different values of \(\lambda_{0}\propto\lambda_{STR}\) using a variable \(d_{tol}\) that initializes \(tol\). \(\mu > 0\) is assigned data-dependently. \(\hat{\mathcal{N}}_{\hat{\xi}}\) is the linear combination of the effective polynomial candidates chosen by \(\mathcal{E}\). \(\mathbb{N}_{\norm{\xi^{STR}}_{0}}\) denotes \(\set{0, 1, \dots, \norm{\xi^{STR}}_{0}-1}\). \(E\) is an indexed set, and \(\boldsymbol{e}_{j}\) is an elementary column vector whose entries are all zero except for the $j^{\text{th}}$ nonzero polynomial candidate. The matrix \(\mathcal{E}\) reduces the dimensionality such that we focus solely on the effective candidates. \(\hat{\xi}\) successively stores the nonzero coefficients in \(\xi^{STR}\). If the library is overcomplete, there exists \(\mathcal{E}\) such that \(\Theta^{\mathcal{M}} \approx P_{k}(\Phi^{\mathcal{M}}(\hat{\theta}))\mathcal{E}\).

The pair values of \((\lambda_1, \lambda_{STR})\) are grid searched with Bayesian information criteria (BIC) \cite{schwarz1978estimating} as the guidance score. The pairs whose PDEs are in agreement with the corresponding preselectors, according to \textbf{Definition \ref{def:agreement}}, are expected.

\begin{definition}[Agreement] \label{def:agreement}
If \(P_{k}\) is regarded as the candidate building function and every nonzero \(f^{\text{th}}_{i+1}\) term can be written as a polynomial of certain \(j^{\text{th}}\) candidates whose \(j^{\text{th}}\) is taken from the set of threshold-passing basis candidate indices \(\set{j \mid I_{j} > \frac{1}{C}}\) (see \textbf{Algorithm \ref{alg:1}}), we determine that the initial discovered PDE of a particular pair of \((\lambda_1, \lambda_{STR})\) is in the ``agreement'' with the \(\lambda_1\)-trained preselector network.
\end{definition}

\noindent The likely models, from which we can voluntarily choose one as the initial discovered PDE, are conceived to be in their agreements and relatively sparse (small \(\norm{\xi^{STR}}_{0}=\abs{E}\)) while conveying sufficiently low BIC scores defined as follows:
\begin{multline} \label{eq:bic}
BIC(\xi^{STR}, \hat{\theta}) = \norm{\xi^{STR}}_{0}\log N_{\mathcal{M}}-2\log \hat{L}(\xi^{STR}, \hat{\theta});\\
\begin{aligned}
\log \hat{L}(\xi^{STR}, \hat{\theta}) = \frac{-N_{\mathcal{M}}}{2}\Bigg(1&+\log 2\pi\\
&+\log \frac{RSS(\xi^{STR}, \hat{\theta})}{N_{\mathcal{M}}}\Bigg),
\end{aligned}\\
RSS(\xi^{STR}, \hat{\theta}) = \sum^{N_{\mathcal{M}}}_{i=1}\abs{\frac{\partial \mathcal{F}_{\hat{\theta}}}{\partial t^{\mathcal{M}}_{i}}-P_{k}(\Phi^{\mathcal{M}}_{i}(\hat{\theta}))\xi^{STR}}^{2}.
\end{multline}
\noindent \(\log \hat{L}(\xi^{STR}, \hat{\theta})\) is the maximized (natural) log-likelihood of the \(\hat{\theta}\)-produced model parameterized by \(\xi^{STR}\). \(RSS\) denotes the real-valued residual sum of squares because the absolute value of each (complex-valued) residual term is considered. BIC formulation is primarily by Statsmodels \cite{seabold2010statsmodels}. The pseudocode for \ref{derivative_preparation} and \ref{initial_pde_identification} is detailed in \textbf{Algorithm \ref{alg:1}}.

Pedagogically, suppose that the preferred initial PDE exemplifies Burgers' PDE; we write the effective candidate matrix concerning the training set of labeled subsamples \(\mathcal{D}\) as
\begin{equation}
\Phi^{\mathcal{D}}_{\mathcal{E}}(\hat{\theta}) = P_{k}(\Phi^{\mathcal{D}}(\hat{\theta}))\mathcal{E} = \mqty[\frac{\partial^{2}\mathcal{F}_{\hat{\theta}}}{\partial {x}^{2}} & \mathcal{F}_{\hat{\theta}}({x}, {t})\frac{\partial\mathcal{F}_{\hat{\theta}}}{\partial {x}}].
\end{equation}

\subsection{dPINNs: Denoising and Finetuning using PINN} \label{dPINNs}
As illustrated by \circled{\textbf{(3)}} of Fig. \ref{fig:overview_framework}, we introduce the denoising PINNs (dPINNs) for achieving the precise recovery of PDE coefficients \(\xi^{*}\) under uncertainties. After \textbf{Algorithm \ref{alg:1}} is performed, we take the weakly physics-constrained solver \(\mathcal{F}_{\hat{\theta}}\) and the initial PDE \(\hat{\mathcal{N}}_{\hat{\xi}}\) to build the dPINNs, minimizing the vigorous physics-informed loss \(\mathcal{L}^{\Tilde{\mathcal{D}}}_{sup}(\hat{\theta})+\mathcal{L}^{\Tilde{\mathcal{D}}^{\prime}}_{unsup}(\hat{\theta}, \hat{\mathcal{N}}_{\hat{\xi}})\) on the denoised dataset \(\Tilde{\mathcal{D}}=\set{(\Tilde{x}_{i}, \Tilde{t}_{i}, \Tilde{u}_{i})^{N_f}_{i=1}}\). The physics loss is generally given by
\begin{equation}
\mathcal{L}^{\Tilde{\mathcal{D}}^{\prime}}_{unsup}(\hat{\theta}, \hat{\mathcal{N}}_{\hat{\xi}}) = \frac{1}{N_f}\sum^{N_f}_{i=1}(\frac{\partial\mathcal{F}_{\hat{\theta}}}{\partial\Tilde{t}_{i}}-\hat{\mathcal{N}}_{\hat{\xi}}[(\Phi^{\Tilde{\mathcal{D}}^{\prime}}_{\mathcal{E}}(\hat{\theta}))_{i}])^{2},
\end{equation}
\noindent where the unsupervised set \(\Tilde{\mathcal{D}}^{\prime}=\set{(\Tilde{x}_{i}, \Tilde{t}_{i})^{N_f}_{i=1}}\) is viewed simply as the slice of \(\Tilde{\mathcal{D}}\) without the supervision. Let us now continue the Burgers' example, we can derive the physics-constraint as 
\begin{equation} \label{eq:burgers_example}
\begin{aligned}
\hat{\mathcal{N}}_{\hat{\xi}}[(\Phi^{\Tilde{\mathcal{D}}^{\prime}}_{\mathcal{E}}(\hat{\theta}))_{i}] &= P_{k}(\Phi^{\Tilde{\mathcal{D}}^{\prime}}_{i}(\hat{\theta}))\mathcal{E}\hat{\xi}\\
&= \hat{\xi}_{1}\frac{\partial^{2}\mathcal{F}_{\hat{\theta}}}{\partial \Tilde{x}^{2}_{i}}+\hat{\xi}_{2}\mathcal{F}_{\hat{\theta}}(\Tilde{x}_i, \Tilde{t}_i)\frac{\partial\mathcal{F}_{\hat{\theta}}}{\partial \Tilde{x}_{i}}.
\end{aligned}
\end{equation}
To continually denoise \(\mathcal{D}\) during the dPINNs' learning, we subtract the transformed noises, initially precomputed by the Discrete Fourier Transform (DFT) algorithm, from both \((x, t)\) and \(u\). The denoising mechanism is formulated as the double affine transformations of the entire training dataset given by
\begin{equation}
\begin{aligned}
(\Tilde{x}, \Tilde{t}) &= (x, t) - \beta_{(x, t)}\odot\mathcal{P}_{\Omega_{(x, t)}}(S_{(x, t)}); \quad S_{(x, t)} = (S_{x}, S_{t}),\\
\Tilde{u} &= u - \beta_{u}\odot\mathcal{P}_{\Omega_{u}}(S_{u}),\\
\end{aligned}
\end{equation}
\noindent where \(\mathcal{P}_{\Omega_{(x, t)}}\) and \(\mathcal{P}_{\Omega_{u}}\) are the projecting functions parameterized by \(\Omega_{(x, t)}\) and \(\Omega_{u}\), capturing the unknown noise distributions. \(\beta_{(x, t)}\) and \(\beta_{u}\) are updated proportional to the unbiased standard deviations \((\sqrt{\V(x)}, \sqrt{\V(t)})\) and \(\sqrt{\V(u)}\), controlling the relevant comparable intensity of the noise corrections. The denoising DFT algorithm, which considers power spectrum density (PSD), is meant to deduct small power frequencies components. The starting noises \(S_{u}\) and \(S_{(x, t)}\) are obtained by limiting frequencies whose power is less than the threshold \(\zeta\). To attain the low-PSD noise for the signal \(\psi \in \set{x, t, u}\), we compute the following quantities:
\begin{multline} \label{eq:denoising_dft}
\begin{aligned}
S_{\psi} &= \psi - \inv{DFT}(DFT^{\zeta}(\psi));\\
DFT^{\zeta}_{k}(\psi) &= \begin{cases}
DFT_{k}(\psi); \quad \textrm{if} \quad PSD_{k}(\psi) > \zeta\\
0; \qquad \qquad \, \, \, \textrm{otherwise},
\end{cases}\\
PSD_{k}(\psi) &= \frac{1}{N_{f}}\norm{DFT_{k}(\psi)}^{2}_{2},\\
\widetilde{PSD}_{k}(\psi) &= \frac{PSD_{k}(\psi)-\E(PSD(\psi))}{\sqrt{\V(PSD(\psi))}},
\end{aligned}\\
\zeta = \E(PSD(\psi))+\alpha\max_{k}(\widetilde{PSD}_{k}(\psi))\sqrt{\V(PSD(\psi))}.
\end{multline}
\noindent Here, \(k\) denotes an index in the frequency domain. \(\zeta\) is defined according to the \(\alpha\) portion of the maximal normalized PSD. \(\E\) and \(\V\) calculates the sample mean and variance over \(k\). We precompute \(S_{({x}, {t})}\) and \(S_{{u}}\), since the gradients cannot flow to \(\alpha\). The denoising physics-informed learning is described in \textbf{Algorithm \ref{alg:2}}. Succeeding the first optimization loop, to compensate the numerical error, least squares (LS) regression (see line \ref{alg:2:ls}) is repeatedly employed on the denoised dataset \(\Tilde{\mathcal{D}}^{\prime}\) until the convergence, i.e., no changes of the optimal unbiased \(\xi^{*}\) are detected between the learning epochs.

\begin{algorithm}
\caption{Denoising physics-informed neural networks' (dPINNs) learning}
\label{alg:2}
\begin{algorithmic}[1]
\State \textbf{Goal:} To achieve the optimal solver parameters $\theta^{*}$ and PDE coefficients $\xi^{*}$.
\State \textbf{Require\footnotemark[1]:} $(x, t)$, $u$, $\hat{\theta}$, $\hat{\mathcal{N}}_{\hat{\xi}}$, initialized $\Omega_{(x, t)}$, $\beta^{\prime}_{(x, t)}$, $\Omega_{u}$ and $\beta^{\prime}_{u}$
\State Compute $S_{({x}, {t})}$, and $S_{{u}}$ using denoising DFT (\ref{eq:denoising_dft})
\State Assign $\beta_{(x, t)}\leftarrow(\sqrt{\V(x)}\beta^{\prime}_{(x, t)}, \sqrt{\V(t)}\beta^{\prime}_{(x, t)})$ \Comment{row vec.}
\State Assign $\beta_{u}\leftarrow\sqrt{\V(u)}\beta^{\prime}_{u}$ \Comment{single parameter}
\While{not converge}\label{alg:2:while_do}
\State Denoise $(\Tilde{x}, \Tilde{t}) \leftarrow (x, t) - \beta_{(x, t)}\odot\mathcal{P}_{\Omega_{(x, t)}}(S_{({x}, {t})})$
\State Denoise $\Tilde{u} \leftarrow u - \beta_{u}\odot\mathcal{P}_{\Omega_{u}}(S_{{u}})$\label{alg:2:denoise_u}
\State Build $\Tilde{\mathcal{D}}^{\prime}\leftarrow\set{(\Tilde{x}_{i}, \Tilde{t}_{i})^{N_f}_{i=1}}$ and $\Tilde{\mathcal{D}}\leftarrow\set{(\Tilde{x}_{i}, \Tilde{t}_{i}, \Tilde{u}_{i})^{N_f}_{i=1}}$\label{alg:2:build_denoised_dataset}
\State Compute loss $\mathcal{L}^{\Tilde{\mathcal{D}}}_{sup}(\hat{\theta})+\mathcal{L}^{\Tilde{\mathcal{D}}^{\prime}}_{unsup}(\hat{\theta}, \hat{\mathcal{N}}_{\hat{\xi}})$ on $\Tilde{\mathcal{D}}$ and $\Tilde{\mathcal{D}}^{\prime}$\label{alg:2:compute_loss}
\State Gradient-based update $\hat{\theta}$, $\hat{\xi}$, $\Omega_{(x, t)}$, $\beta^{\prime}_{(x, t)}$, $\Omega_{u}$ and $\beta^{\prime}_{u}$\label{alg:2:grad_update}
\EndWhile\label{alg:2:end_while}
\State Minimize $\mathcal{L}^{\Tilde{\mathcal{D}}}_{sup}(\hat{\theta})+\mathcal{L}^{\Tilde{\mathcal{D}}^{\prime}}_{unsup}(\hat{\theta}, \hat{\mathcal{N}}_{\xi^{*}})$; $\hat{\mathcal{N}}_{\xi^{*}}$ is represented by $\xi^{*}\leftarrow\inv{((\Phi^{\Tilde{\mathcal{D}}^{\prime}}_{\mathcal{E}}(\hat{\theta}))^{\intercal}\Phi^{\Tilde{\mathcal{D}}^{\prime}}_{\mathcal{E}}(\hat{\theta}))}(\Phi^{\Tilde{\mathcal{D}}^{\prime}}_{\mathcal{E}}(\hat{\theta}))^{\intercal}\frac{\partial\mathcal{F}_{\hat{\theta}}}{\partial \Tilde{t}}$\Comment{Redo line \ref{alg:2:while_do}-\ref{alg:2:end_while} with $\xi^{*}$ iteratively resolved between line \ref{alg:2:build_denoised_dataset} and \ref{alg:2:compute_loss} by LS instead of its gradient-based update at line \ref{alg:2:grad_update}.}\label{alg:2:ls}
\State \textbf{Return\footnotemark[2]:} $(x^{*}, t^{*})$, $u^{*}$, $\theta^{*}$, $\xi^{*}$, $\Omega^{*}_{(x, t)}$, $\beta^{*}_{(x, t)}$, $\Omega^{*}_{u}$ and $\beta^{*}_{u}$
\end{algorithmic}
\hrulefill \\
\footnotesize \footnotemark[1]{$\hat{\theta}$ and $\hat{\mathcal{N}}_{\hat{\xi}}$ are attained from \textbf{Algorithm \ref{alg:1}.}}
\footnotesize \footnotemark[2]{The learned outputs are assigned as the optimal parameters superscripted with the asterisk ($*$) notation.}
\end{algorithm}
\clearpage
\section{Experiments and Results} \label{results}
We experimented with 5 canonical PDEs, including 3 ordinary PDEs and 2 complex-valued PDEs, to investigate the accuracy and robustness of our proposed method. We present the results of \circled{\textbf{(1)}} Derivative preparation and \circled{\textbf{(2)}} Initial PDE discovery and discuss the regularization hyperparameter effects on finding the appropriate initial PDE. Later, we show the tolerance of \circled{\textbf{(3)}} dPINNs against noise in both \((x, t)\)\&\(u\) for each PDE as well as against the decreasing number of training samples (scarce data). Beyond the numerical results, we visualize how the projection networks handle the increasing noise intensity in the exemplar of discovering Burgers' PDE.

\subsection{Canonical PDEs}
\subsubsection{Burgers' PDE}
The equation arises in various areas of applied mathematics such as fluid mechanics and traffic flow \cite{basdevant1986spectral}. We consider the following Burgers' equation dataset simulated with Dirichlet boundary conditions, studied in \cite{raissi2019physics}.
\begin{equation} \label{eq:burgers_pde}
\begin{gathered}
    u_t + uu_x - \nu u_{xx} = 0; \quad \nu = \frac{0.01}{\pi},\, x \in [-1, 1],\, t \in [0, 1].
\end{gathered}
\end{equation}
\noindent Different from the previous works such as \cite{rudy2017data, raissi2018deep} where the viscosity of fluid \(\nu\), was set to \(0.1\); thus, the smooth fluid speed without a shock wave, here \(\nu\) = \(\frac{0.01}{\pi}\) is so small that the shock wave emerges.

\subsubsection{Korteweg–De Vries (KdV) PDE}
The KdV equation \cite{korteweg1895xli} is a nonlinear dispersive PDE for describing the motion of unidirectional shallow water surfaces. For a function \(u(x, t)\) the actual form of KdV we consider is expressed as
\begin{equation}
\begin{gathered}
    u_t + 6uu_x + u_{xxx} = 0; \quad x \in [0, 50],\, t \in [0, 50].
\end{gathered}
\end{equation}
\noindent KdV was known to have soliton solutions, representing two one-way moving waves with different amplitudes. Such characteristics challenge discovery methods to distinguish and yield the sparsest governing PDE that generalizes the situation. The PDE is also an excellent prototypical example to test discovering the relatively high-order spatial derivative \(u_{xxx}\).

\subsubsection{Kuramoto–Sivashinsky (KS) PDE}
The KS or flame equation is a chaotic nonlinear PDE with a spatial fourth-order derivative term, primarily to model the diffusive instabilities in a laminar flow. The PDE reads
\begin{equation}
\begin{gathered}
    u_t + uu_x + u_{xx} + u_{xxxx} = 0; \quad x \in [0, 100],\, t \in [0, 100].
\end{gathered}
\end{equation}
\noindent The solution was generated with an initial condition \(u(x, 0 ) = \cos(\frac{x}{16})(1+\sin(\frac{x}{16}))\), integrated up to the wide temporal bound of \([0, 100]\) \cite{rudy2017data}. Consequently, we got a chaotic and complicated PDE solution. Raissi \cite{raissi2018deep} very first noticed that it was challenging to fit a vanilla neural network to the entire chaotic solution while minimizing the residual physics loss; for example, \(\min_{\theta, \theta_{s}}(\mathcal{L}^{\mathcal{D}}_{sup}(\theta)+\mathcal{L}^{\mathcal{D}_{s}}_{unsup}(\theta, \theta_{s}))\). A similar problem was independently found by Rudy {\it et al}. \cite{rudy2017data} that when encountering the whole chaotic domain of KS, the PDEs produced by STRidge could be inaccurate and unstable with the complication of noise.

\subsubsection{Quantum Harmonic Oscillator (QHO) PDE}
The quantum harmonic oscillator is the Schrodinger equation with a parabolic potential \(0.5x^{2}\). The PDE is given by
\begin{equation}
\begin{gathered}
    iu_t + \frac{1}{2}u_{xx} - \frac{x^{2}}{2}u = 0; \quad x \in [-7.5, 7.5],\, t \in [0, 4].
\end{gathered}
\end{equation}
\noindent Following \cite{rudy2017data}, we construct the basis candidate matrix that includes the parabolic potential.

\subsubsection{Nonlinear Schrodinger (NLS) PDE}
The nonlinear Schrodinger equation is used to study nonlinear wave propagation. The true discretization studied in \cite{raissi2019physics}, is expressed by
\begin{equation}
\begin{gathered}
    iu_t + \frac{1}{2}u_{xx} + u\norm{u}^{2}_{2} = 0; \quad x \in [-5, 5],\, t \in [0, \frac{\pi}{2}].
\end{gathered}
\end{equation}
\noindent We include candidate terms depending on the magnitude of the solution, e.g., \(\norm{u}^{2}_{2}\), which may appear in the correct identification of the dynamics of the complex-valued function.

\subsection{Experimental Settings} \label{exp_settings}
The training data points \((x, t)\)\&\(u\) are randomly subsampled from all the generated discretized points in the domain according to the size \(N_f\) specified in Table \ref{tab:noise}. All the discretized (noisy) data points are exploited as the validation set for early stopping once the validation MSE drops during pretraining and converging the solver network that minimizes the MSE loss. \(N_r=(1, 1, 0.5, 0.5, 1)N_f\) for Burgers', KdV, KS, QHO and NLS PDE, respectively. The solver architecture comprises 6 hidden layers with 50 neurons each and Tanh activation functions in the between. For the preselector, \(W^b\) are devised as a single hidden layer. At the same time, the rest parameters \(\theta^{r}_s\) are implemented as a sequence of 3 hidden layers, each with 50 neurons whose outputs are layer normalized \cite{ba2016layer}, randomly dropped out \cite{srivastava2014dropout} and Tanh activated, excluding Tanh from the final layer. The dropout probability is \(0.1\) for KdV and KS, otherwise is 0.0. Hidden weights are initialized by uniform Xavier \cite{glorot2010understanding} and biases are initialized to \(0.01\). \(\sigma(\cdot) = \frac{1}{2}(\tanh{(\cdot)}+1)\) is defined for all the canonical models except for Burgers' PDE, \(\sigma(\cdot)=\frac{1}{1 + \exp(-1(\cdot))}\), Sigmoid is employed to convey the flexibility in the design. \(\lambda_1\) is varied for accomplishing the suitable value while \(\lambda_2\) is set to \(0.1\). The projection networks \(\Omega_{(x, t)}\) and \(\Omega_{u}\) are 2 hidden layers, each having 32 neurons with Tanh; hence, the final layer's raw outputs of the networks \(\mathcal{P}_{\Omega_{(x, t)}}\) and \(\mathcal{P}_{\Omega_{u}}\) are activated by Tanh. \(\beta^{\prime}_{(x, t)}\) and \(\beta^{\prime}_{u}\) are initialized at \(10^{-3}\) for the ordinary PDEs and \(10^{-5}\) for the complex-valued PDEs (QHO and NLS).

For \textbf{Algorithm \ref{alg:1}}, full-batch stochastic LBFGS \cite{yatawatta2019stochastic} and vanilla LBFGS \cite{liu1989limited}, with 0.1 step sizes and the strong Wolfe line search, are leveraged separately, to pretrain and converge the solver network. The pretraining (second-order optimization) epoch is limited to 1 to prevent overfitting in the noisy \((x, t)\)\&\(u\) case. MADGRAD \cite{defazio2021adaptivity} with gradient-deconflicting PCGrad \cite{yu2020gradient} is applied to joint learn (line \ref{alg:1:joint_train}) for 1,000 epochs in Burgers' and KdV cases. The weighted average with the ratios \(\mathcal{L}^{\mathcal{D}}_{sup}(\theta):(\mathcal{L}^{\mathcal{D}_{s}}_{unsup}(\theta, \theta_{s})+R^{\mathcal{D}_{s}}_{\eta}(\theta, W^b)) = 1:1\) and \(1:10^{-3}\) are put to optimize for 300 and 1,500 epochs in KS and the complex-valued PDEs. The learning rate for updating the pretrained \(\theta\) is assigned with a low value of \(10^{-7}\), while the higher rates from \((10^{-2}, 10^{-2}, 10^{-3}, 10^{-1}, 10^{-1})\) are set for updating untrained \(\theta_s\). \(\kappa\) is set, in the same dataset order, to \((0.75, 0.7, 0.8, 0.9, 0.9)\) before the first gradient updates of the joint training. Then, LBFGS \cite{liu1989limited} is mainly used for the dPINNs' learning (\textbf{Algorithm \ref{alg:2}}). For every noisy KdV and KS experimental case, the denoising-related parameters $\Omega_{(x, t)}$, $\beta^{\prime}_{(x, t)}$, $\Omega_{u}$ and $\beta^{\prime}_{u}$ are reinitialized with the conceivably closer estimate \(\hat{\theta}\) prior to executing the subroutine at line \ref{alg:2:ls}.

As for the input of STRidge, the candidate library is \(P_{2}(\cdot)\), collecting unbiased interaction-only real-valued polynomial features up to the \nth{2} degree of the estimated PDE solution and its partial derivatives \(P_{2}(\cdot)\), computed with respect to \(\mathcal{M}\).

The precomputed denoising DFT is configured with \(\alpha=0.1\) for all the canonical PDEs. \(DFT\) and \(\inv{DFT}\) (the inverse transform) are the one-dimensional fft and ifft operators from PyTorch \cite{paszke2019pytorch} package. Our nPIML framework is as well implemented dominantly using PyTorch package.

In the noisy experiments, we presume that a matrix, say \(z\), gets perturbed, right after the time of its subsampling, by the p\% biased (no Bessel's correction) standard deviation (\(std\)) of Gaussian noise \(Z\) simulated as follows:
\begin{equation}
noise(z, p) = \frac{p\cdot std(z)}{100} \times Z;\, \forall i,j(Z_{ij} \sim \boldsymbol{\mathcal{N}}(0, 1)).
\end{equation}
\noindent Suppose that 1\% noise is exerted, subsampled \(u\) and \((x, t)\) get polluted in turn with \(noise(u, 1)\) and \((\frac{noise(x, 1)}{\sqrt{2}}, \frac{noise(t, 1)}{\sqrt{2}})\).

The metric to measure how far an estimate \(\xi^{est}\) from the ground truth \(\xi\) is \(mean(\delta) \pm std(\delta)\) over all \(j\) effective coefficients in \(\xi^{est}\). If only the correct candidates are identified, \(\delta_{j}=\delta_{j}(\xi^{est}, \xi)\) is the \%coefficient error (\%CE) defined as 
\begin{equation} \label{eq:percent_ce}
\delta_{j} = \abs{\xi^{est}_{j}-\xi_{j}}/\abs{\xi_{j}} \times 100\%;\, j \in \set{1, \dots, cols(\Theta)}.
\end{equation}
\noindent In Table \ref{tab:noise}, \ref{tab:tolerance} and \ref{tab:high_noise}, \(\xi^{est} \in \set{\hat{\xi}, \xi^{*}}\). \(cols(\Theta)\) represents the number of column(s) of \(\Theta\).

\textit{Specific Treatments for Complex-valued PDEs:} 
Our complex neural networks are initialized based on the prior work called Deep complex networks \cite{trabelsi2018deep}. Since the spatio-temporal points lay on a real 2-dimensional plane, the model starts from 1 (real) hidden layer with 200 neurons, followed by 5 complex linear layers, each consisting of 200 neurons that account for 100 real parameters and 100 imaginary parameters. Note that the complex forward pass is essentially iteratively performing naive complex-valued matrix multiplication and bias addition. The differentiation of complex-valued \(\mathcal{F}_{\theta}(x, t)\), respecting a real-valued vector, e.g., \(x\), can be computed distributively. Concretely, we apply automatic differentiation to the real and imaginary parts with respect to \(x\) separately; then, we form the output complex-valued matrix as 
\begin{equation}
\frac{\partial \mathcal{F}_{\theta}(x, t)}{\partial x} = \frac{\partial \Re(\mathcal{F}_{\theta}(x, t))}{\partial x} + \frac{\partial \Im(\mathcal{F}_{\theta}(x, t))}{\partial x}i;\, i^{2} = -1.
\end{equation}
Likewise, \(W^b\) of the preselector is treated as a single complex linear layer, including the bias, with 50 neurons. \(\theta^{r}_s\) is modeled by 3 complex linear layers, each with total 50 neurons that are batch normalized \cite{ioffe2015batch} and component-wise Relu activated.

Because the estimated PDE solution is in complex form, we may include norm-based atomic candidates, e.g., \(\norm{\mathcal{F}_{\hat{\theta}}(x^{\mathcal{M}}, t^{\mathcal{M}})}^{2}_{2}\), on which the all (not interaction-only) polynomial features, up to the \nth{2} degree, are built. Once prepared, the candidate library can be directly input to STRidge.

\subsection{Effect of Regularization Hyperparameters on Initial PDE Identification}
For each canonical PDE, we present the domain of interest from which the metadata \(\mathcal{M}\) is generated for the initial PDE extraction. We then concentrate on the multi-perspective assessment of the different discovered PDEs by STRidge while varying the two major regularization hyperparameters: \(\lambda_{1}\) of the preselector network and \(\lambda_{STR}\) of STRidge algorithm. Before the finetuning process, we present how accurate the initial discovered PDEs in order, concerning the following three cases distinguished by the noise conditions: noiseless dataset, noiseless \((x, t)\) but noisy \(u\), and noisy \((x, t)\)\&\(u\) in which the spatial-temporal \((x, t)\) becomes mesh-free.

\begin{figure}
\centering
\includegraphics[width=0.5\textwidth]{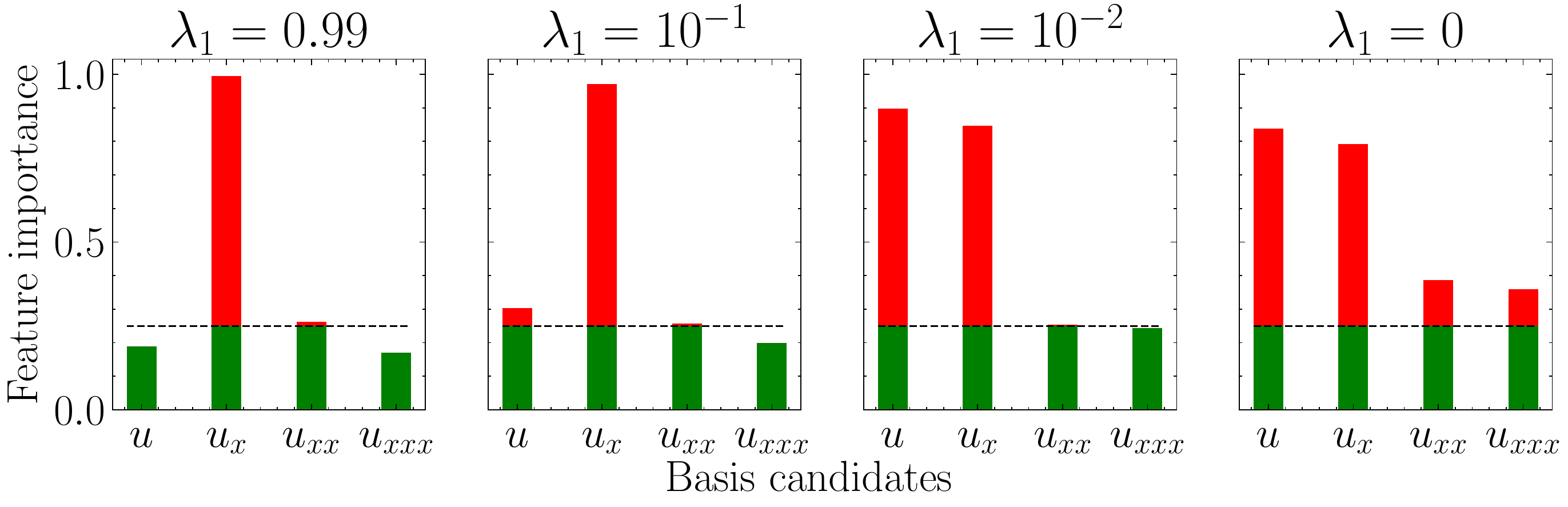}
\caption{Burgers: Learned feature importance with varied \(\lambda_{1}\)}
\vspace*{-\baselineskip}
\label{fig:feature_importance_burgers}
\end{figure}

\subsubsection{Initial Discovered Burgers' PDE} \label{discovering_burgers}
We trained the preselector network with varying \(\lambda_1\) to perceive the significance of each candidate. The distributed feature importance values (\(I_j\) for each \(j^{\text{th}}\) basis candidate) are presented in Fig. \ref{fig:feature_importance_burgers}. Although several choices of the expressive subset of passing-threshold candidates are contributed, identifying the optimal set is still not obvious by merely adjusting \(\lambda_1\). Hence, STRidge was subsequently employed multiple times with diverse levels of regularization intensity \(\lambda_{STR}\). For convenience, we simply set \(\forall i \leqslant N_f+N_r, (x^{\mathcal{M}}_{i}, t^{\mathcal{M}}_{i}) = (x_{i}, t_{i})\) for all Burgers' experimental cases that differed in the noise conditions. The cross results, Table \ref{tab:lambdas}, are assessed for obtaining the initial discovered PDE that is preferably conceived of being agreed with the corresponding preselector and sparse with a sufficiently low BIC score. We could have imposed an explicit metric for selecting the best initial governing PDE, but we did not due to the no-free-lunch problem of defining the single criterion that always determines the actual function of every physical system; therefore, the optimality subject to one's wilfulness.

\begin{table}
\begin{center}
\begin{tabular}{cccc}
\hline\noalign{\smallskip}
$\lambda_{1}$/$\lambda_{STR}$ & $10^{-6}$ & $10^{-3}$ & $10^{0}$\\
\noalign{\smallskip}\hline\noalign{\smallskip}
$0.99$ & \makecell[l]{$[u_{xx}, uu_x,$ \\ $\quad uu_{xxx}, u_xu_{xx}]$} & \makecell[c]{$[u_{xx}, uu_x]$} & \makecell[c]{$[uu_x]$}\\
 & (\textbf{-8,723.69}) & (-7,636.39) & (15,823.14)\\
\noalign{\smallskip}\hline\noalign{\smallskip}
$10^{-1}$ & \makecell[l]{$[u_{xx}, uu_x,$ \\ $\quad uu_{xxx}, u_xu_{xx}]$} & \makecell[c]{$[u_{xx}, uu_x]$} & \makecell[c]{$[uu_x]$}\\
 & (\textbf{-8,456.28}) & (\textcolor{blue}{-7,154.65}) \checkmark & (\textcolor{blue}{15,824.98})\\
\noalign{\smallskip}\hline\noalign{\smallskip}
$10^{-2}$ & \makecell[l]{$[u_{xx}, uu_x,$ \\ $\quad uu_{xxx}, u_xu_{xx}]$} & \makecell[c]{$[u_{xx}, uu_x]$} & \makecell[c]{$[uu_x]$}\\
 & (\textbf{-8,294.55}) & (\textcolor{blue}{-7,178.84}) \pmb{\checkmark} & (\textcolor{blue}{15,824.29})\\
\noalign{\smallskip}\hline
\noalign{\smallskip}\hline\noalign{\smallskip}
\makecell[c]{$0$ \\ (Supplement)} & \makecell[l]{$[u_{xx}, uu_x,$ \\ $\quad uu_{xxx}, u_xu_{xx}]$} & \makecell[c]{$[u_{xx}, uu_x]$} & \makecell[c]{$[uu_x]$}\\
 & (\textbf{\textcolor{blue}{-8,437.81}}) \checkmark & (\textcolor{blue}{-7,243.32}) \checkmark & (\textcolor{blue}{15,827.68})\\
\noalign{\smallskip}\hline
\end{tabular}
\caption{\textbf{Burgers regularization hyperparameter selection}: Concerning the coefficient selection criteria, STRidge's \(\lambda_{0}\), controlling the \(L_{0}\)-penalty, is set to \(10^{4}\lambda_{STR}\varepsilon\), and \(d_{tol}\) equals \(2\) for the three noise conditions. The assignment of \((\mu, \lambda_{STR}, d_{tol})\) is purely for gathering the likely different PDEs. Each PDE is accompanied by the ``(BIC)'' score. \textcolor{blue}{Blue} indicates the agreement. \textbf{Bold} means the lowest BIC score, compared to the scores acquired by the same \(\lambda_{1}\). Among the agreed models, we check (\checkmark) the sparse PDEs with \(\norm{\xi^{STR}}_{0} \leqslant 4\), which demonstrate sufficiently low BIC score. The PDE with \pmb{\checkmark} is regarded as the initial guess.}
\vspace*{-\baselineskip}
\vspace*{-\baselineskip}
\label{tab:lambdas}
\end{center}
\end{table}

Assigning the $\lambda_{1} = 0.99$ is so high that the true candidate, i.e., \(u\), is lacking from the passing-threshold candidates. Accordingly, the resulting PDEs cannot match the particular importance scores. The preselector properly focuses on the true candidates when \(\lambda_{1}\) is set to \(10^{-1}\) and \(10^{-2}\). Notice that \(u_{xx}\) consistently passes the threshold with marginal values, conveying the small viscosity estimates. As seen in Table \ref{tab:lambdas}, for \(\lambda_{1}>0\), \(10^{-2}\) gave the best initial result, covering the sparse PDE with the lowest BIC among the agreed models. For a new real-world problem without any knowledge about the underlying equation, we advise selecting a \(\lambda_{1}\) that cuts out some potentially unimportant candidates and causes the agreement with the Pareto-optimal solution suggested by STRidge, e.g., the one that minimizes \(\frac{\Delta BIC}{\Delta \norm{\xi^{STR}}_{0}}\).

Deciding on the value of \(\lambda_{STR}\) requires an akin principle: the values that are too low or high are likely to yield incorrect forms. For example, $\lambda_{STR} = 10^{0}$ is immensely high, outputting the too sparse and noninformative PDE with the single effective \(uu_x\), delivering the high BIC scores. $\lambda_{STR} = 10^{-3}$ is more suitable, suggesting the sparse models, which conform with the preselectors and offer the low BIC scores that vastly improve from those given by $\lambda_{STR} = 10^{0}$. Conditioned by \(\lambda_1>0\), $u_t = 0.003063u_{xx}-0.986174uu_{x}$ contains the few terms and offers the minimal BIC among the acceptable PDEs; thus, taken as our initial guess (\pmb{\checkmark}) to be finetuned. Remind that, when comparing the models from diverse values of \(\lambda_1\), although their functions differ solely in the set of PDE coefficients, they cannot be directly compared because the change in \(\hat{\theta}\) affects \(\mathcal{F}_{\hat{\theta}}(\cdot)\), i.e. \(\frac{\partial \mathcal{F}_{\hat{\theta}}}{\partial t^{\mathcal{M}}}\) varies (see (\ref{eq:bic})); therefore the slightly flustered RSS scales without an explicit static referenced time derivative. Nevertheless, we straightforwardly prefer the one with the lower BIC score. By the disagreements, the sparsity-promoting preselectors trained with \(\lambda_{1}>0\) all entails that $\lambda_{STR} = 10^{-6}$ gives overly parameterized models, with the minor improvements per the increased independent candidates. If we were to independently have the mere consideration on $\lambda_{1} = 0$ or technically diminutive to a certain value, none of the basis candidates would probably get deselected, and the resulted PDEs would be all in their agreements. The justification, whether including \(uu_{xxx}\) and \(u_{x}u_{xx}\) worth the reduction in BIC, would turn ambiguous, though the PDE outcome by \((\lambda_1, \lambda_{STR})=(0, 10^{-3})\): $u_t = 0.003063u_{xx}-0.985882uu_{x}$ captures the ground on par with our PDE guess (\(\pmb{\checkmark}\)). If the preselector were not at all constructed, the concern would still persist. For the noisy cases, the \%CE (see (\ref{eq:percent_ce})) of the initial PDE estimates are listed in the nPIML: IPI row of Table \ref{tab:noise}.

\begin{figure}
\centering
\includegraphics[width=0.5\textwidth]{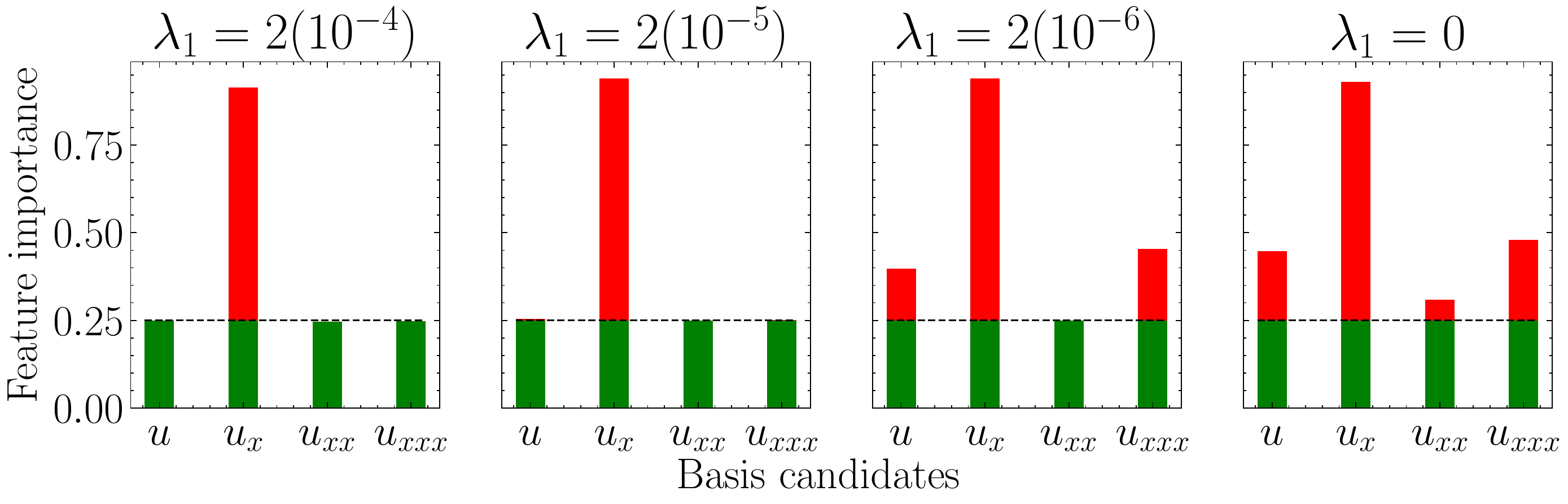}
\caption{KdV: Learned feature importance with varied \(\lambda_{1}\)}
\vspace*{-\baselineskip}
\label{fig:feature_importance_kdv}
\end{figure}

\begin{table}
\begin{center}
\begin{tabular}{cccc}
\hline\noalign{\smallskip}
$\lambda_{1}$/$\lambda_{STR}$ & $10^{-5}$ & $10^{-3}$ & $10^{-1}$\\
\noalign{\smallskip}\hline\noalign{\smallskip}
$2(10^{-4})$ & \makecell[l]{$[u_x, u_{xxx}, uu_x,$ \\ $\quad uu_{xxx}, u_{x}u_{xx}]$} & \makecell[l]{$[u_{xxx}, uu_{x}]$} & \makecell[c]{$[u_x]$}\\
 & (\textbf{-651,496.23}) & (-593,260.84) & (\textcolor{blue}{-493,869.28})\\
\noalign{\smallskip}\hline\noalign{\smallskip}
$2(10^{-5})$ & \makecell[l]{$[u_x, u_{xxx}, uu_x,$ \\ $\quad uu_{xxx}, u_{x}u_{xx}]$} & \makecell[c]{$[u_{xxx}, uu_{x}]$} & \makecell[c]{$[u_x]$}\\
 & (\textbf{-651,650.73}) & (\textcolor{blue}{-593,259.27}) \checkmark & (\textcolor{blue}{-493,885.29})\\
\noalign{\smallskip}\hline\noalign{\smallskip}
$2(10^{-6})$ & \makecell[l]{$[u_x, u_{xxx}, uu_x,$ \\ $\quad uu_{xxx}, u_{x}u_{xx}]$} & \makecell[c]{$[u_{xxx}, uu_{x}]$} & \makecell[c]{$[u_x]$}\\
 & (\textbf{-651,782.07}) & (\textcolor{blue}{-593,389.01}) \pmb{\checkmark} & (\textcolor{blue}{-493,868.73})\\
\noalign{\smallskip}\hline
\noalign{\smallskip}\hline\noalign{\smallskip}
\makecell[c]{$0$ \\ (Supplement)} & \makecell[l]{$[u_x, u_{xxx}, uu_x,$ \\ $\quad uu_{xxx}, u_{x}u_{xx}]$} & \makecell[c]{$[u_{xxx}, uu_{x}]$} & \makecell[c]{$[u_x]$}\\
 & (\textbf{\textcolor{blue}{-651,733.37}}) & (\textcolor{blue}{-593,275.19}) \checkmark & (\textcolor{blue}{-493,851.71})\\
\noalign{\smallskip}\hline
\end{tabular}
\caption{\textbf{KdV regularization hyperparameter selection}: STRidge's \(\lambda_{0}\) is set to \(10^{2}\lambda_{STR}\varepsilon\), and \(d_{tol}\) equals \(1\) for the three noise conditions.}
\vspace*{-\baselineskip}
\vspace*{-\baselineskip}
\label{tab:lambdas_kdv}
\end{center}
\end{table}

\subsubsection{Initial Discovered KdV PDE} \label{discovering_kdv}
We inspect how the preselector weights each basis candidate in Fig. \ref{fig:feature_importance_kdv}. Trained with \(\lambda_{1}=2(10^{-5})\) or \(2(10^{-6})\), the preselector can capture the true candidates while the relatively high value of \(\lambda_{1}=2(10^{-4})\) solely let \(u_x\) pass the threshold. \(u\) and \(u_{xxx}\) barely pass the threshold if \(\lambda_{1}=2(10^{-5})\), nonetheless their effectiveness become vivid when \(\lambda_{1} \leqslant 2(10^{-6})\).

STRidge was leveraged multiple times on the candidate library built on \(\mathcal{M}\). For KdV, we regarded the metadata as the linear discretization of the entire spatio-temporal domain; \(N_{\mathcal{M}}=64,128\), facilitating the disambiguation of the different wave amplitudes. The found PDEs for the several pair of \((\lambda_1, \lambda_{STR})\) are listed in Table \ref{tab:lambdas_kdv}. By pondering the PDEs that harmonize with \(\lambda_1>0\), we neglect the selection of the PDEs with the minimal BIC (for a particular \(\lambda_{1}\)) because they neither agree with the \(L_0\)-penalized feature importance nor be sparse as expected. The reduced BIC per an increasing effective term of transition from \(\lambda_{STR}=10^{-3}\) to \(\lambda_{STR}=10^{-5}\) is much less when compared with moving from \(\lambda_{STR}=10^{-1}\) to \(\lambda_{STR}=10^{-3}\), signifying the inefficiency of including the unnecessary terms. Remark that setting \(\lambda_{STR}=10^{-1}\) gives the PDEs, each describing a one-way traveling wave which can be considered as the relaxed form of KdV PDE, still not well fit the overall character of the dataset. Based on the mentioned justification, we thus prefer \(\lambda_{STR}=10^{-3}\), and choose the agreed PDE with the better BIC, taking the form of \(u_t = -0.989065u_{xxx}-5.961087uu_{x}\) as our initial guess (\pmb{\checkmark}). The selected PDE is noticed as a more precise to the ground truth than the PDE based \(\lambda_1=0\), which is \(u_t = -0.988350u_{xxx}-5.959614uu_{x}\). Also, just naively, the BIC cannot elucidate the overfitting hurdle without the auxiliary knowledge gained by varying \(\lambda_1>0\). For the noisy KdV cases, the initial results \%CE of the \textbf{Algorithm \ref{alg:1}} are as well shown in the nPIML: IPI row of in Table \ref{tab:noise}.

\begin{figure}
\centering
\includegraphics[width=0.45\textwidth]{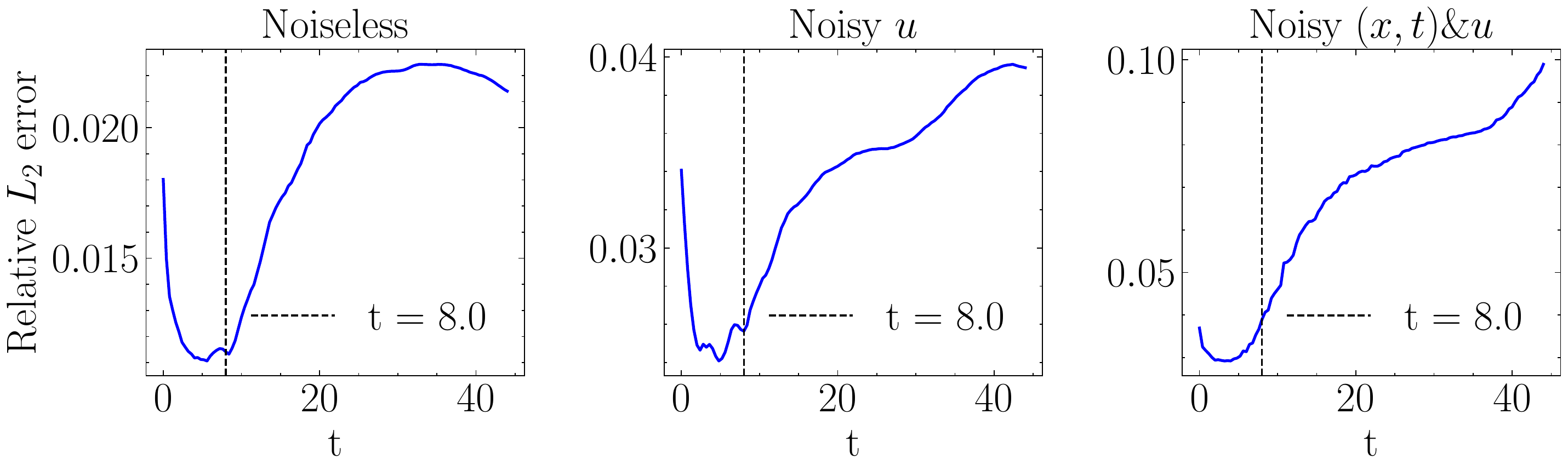}
\caption{KS: Training relative \(L_2\) error of the learned (from \(N_f=80,000\)) solver \(\hat{\theta}\) against temporally varying sub-regions of the KS training set bounded by \([0, 100] \times [0, 44]\), revealing a local optimum around the stability domain at the beginning of the evolution.}
\label{fig:loss_plots_ks}
\end{figure}

\begin{figure}
\centering
\includegraphics[width=0.45\textwidth]{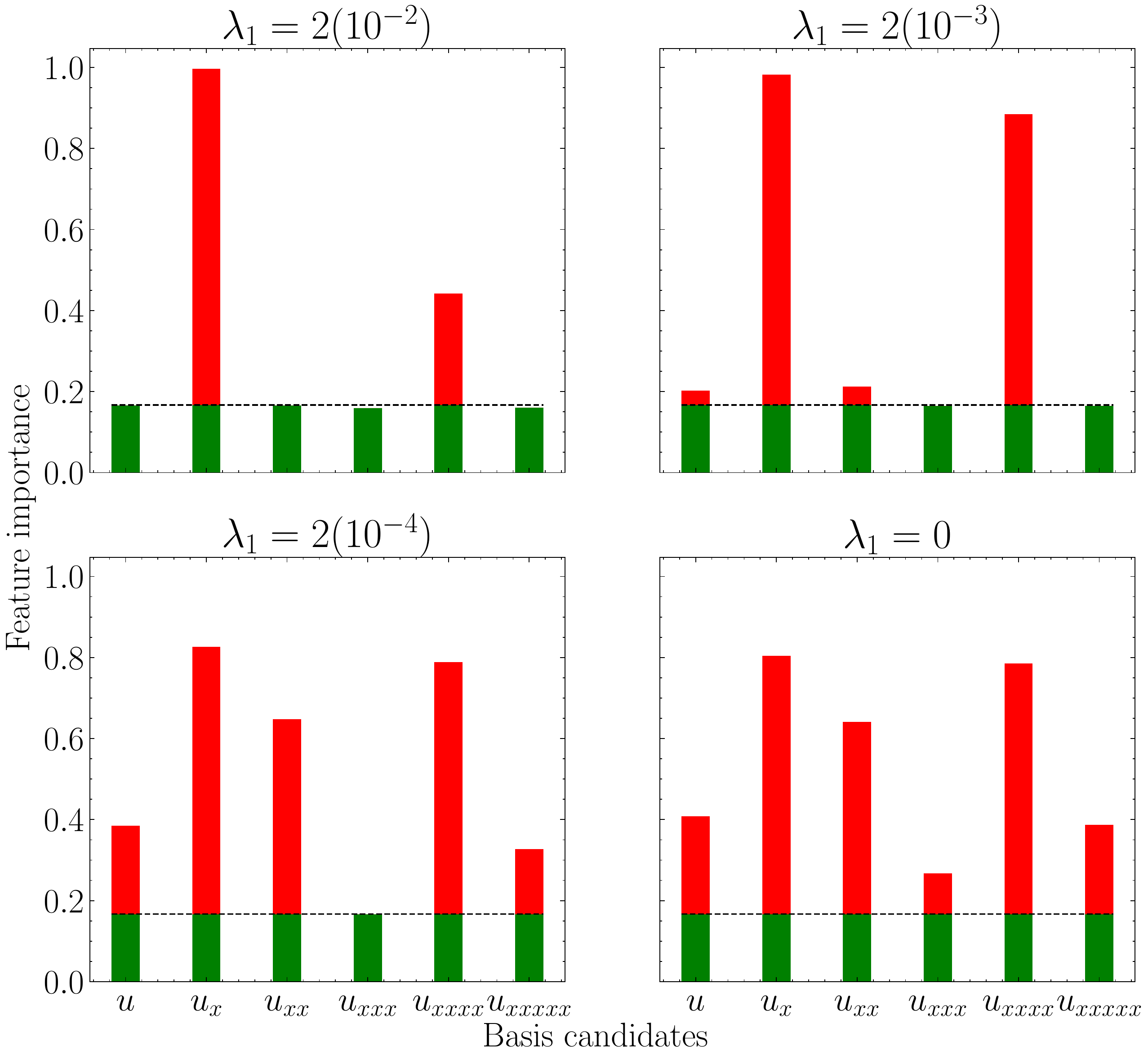}
\caption{KS: Learned feature importance with varied \(\lambda_{1}\)}
\vspace*{-\baselineskip}\label{fig:feature_importance_ks}
\end{figure}

\subsubsection{Initial Discovered KS PDE} \label{discovering_ks}
Our early attempt was performing \textbf{Algorithm \ref{alg:1}} with train/validation sets. The training samples were abundant as \(N_f=80,000\). \(N_r=0\) was chosen to avert the overflow of 48,601 MiB GPU memory because of the computation up to the fifth-order \(u_{xxxxx}\). Unfortunately, suggested by the plots in Fig. \ref{fig:loss_plots_ks}, we have quickly realized that the relative \(L_2\) error of the solver network starts diverging, especially if noise exists when entering the highly chaotic region of KS, admonishing the evidence of training PINN burdensome upon the full-field domain \cite{raissi2018deep}. The issue leads to unreliable derivation estimation; hence, the non-sparse and cluttered discoveries of the governing function by STRidge.

We bypass the complication by selectively focusing on the samples from a more stable sub-region at the beginning of the evolution, where the solver can accurately approximate as indicated by the relative \(L_2\) error plots in Fig. \ref{fig:loss_plots_ks}. We assumed that the unknown PDE governs persistently throughout the evolution; nevertheless, the presumption does not universally hold since specific coefficients of the chaotic behavior can be distinct over time \cite{quade2018sparse}. Based on the encountered evidences, as a result, the first 21,504 (\(1,024\times21\)) discretized points within \([0, 100]\times[0,8]\), were instead used with randomly generated nonoverlapping 10,752 unsupervised points for the (re)training in the noiseless experiment. The temporally-wise increased number of training samples to be the first 30,000 polluted discretized points, where \(t \leqslant 11.6\), were used with randomly generated disjoint 15,000 unsupervised points for both the noisy experiments. The validation sets were homogeneously left unaffected. Before the initial PDE identification, we retrained the networks using \textbf{Algorithm \ref{alg:1}} once from scratch on these altered, better stability training sets.

We investigate the learned feature importance of the preselector for ranking each potential atomic candidate, helping us choose the right PDE as presented in Fig. \ref{fig:feature_importance_ks}. It is intriguing to discern that \(u_{xxxx}\) is one of the essential terms for every choice of \(\lambda_{1}\), despite its order being 4, implying that the high-order derivative is plausible to be included.

\begin{table}[t]
\begin{center}
\addtolength{\tabcolsep}{-3pt}
\begin{tabular}{cccc}
\hline\noalign{\smallskip}
$\lambda_{1}$/$\lambda_{STR}$ & $10^{-5}$ & $10^{-3}$ & $10^{-1}$\\
\noalign{\smallskip}\hline\noalign{\smallskip}
$2(10^{-2})$ & \makecell[c]{$[u_{xx}, u_{xxxx}, uu_x,$ \\ $uu_{xxx}, uu_{xxxxx},$ \\ $u_{x}u_{xx}, u_{xx}u_{xxx},$ \\ $u_{xx}u_{xxxxx}]$} & \makecell[c]{$[u_{xx}, u_{xxxx}, uu_x]$} & \makecell[c]{$[uu_x]$}\\
 & (\textbf{-153,326.24}) & (-141,117.36) & (-67,989.30)\\
\noalign{\smallskip}\hline\noalign{\smallskip}
$2(10^{-3})$ & \makecell[c]{$[u_{xx}, u_{xxxx}, uu_x,$ \\ $uu_{xxx}, uu_{xxxxx},$ \\ $u_{x}u_{xx}, u_{xx}u_{xxx},$ \\ $u_{xx}u_{xxxxx}]$} & \makecell[c]{$[u_{xx}, u_{xxxx}, uu_x]$} & \makecell[c]{$[uu_x]$}\\
 & (\textbf{-153,661.00}) & (\textcolor{blue}{-141,032.21}) \pmb{\checkmark} & (\textcolor{blue}{-67,956.02})\\
\noalign{\smallskip}\hline\noalign{\smallskip}
$2(10^{-4})$ & \makecell[c]{$[u_{xx}, u_{xxxx}, uu_x,$ \\ $uu_{xxx}, uu_{xxxxx},$ \\ $u_{x}u_{xx}, u_{xx}u_{xxx},$ \\ $u_{xx}u_{xxxxx}]$} & \makecell[c]{$[u_{xx}, u_{xxxx}, uu_x]$} & \makecell[c]{$[uu_x]$}\\
 & (\textbf{-151,328.93}) & (\textcolor{blue}{-138,842.33}) \checkmark & (\textcolor{blue}{-68,022.47})\\
\noalign{\smallskip}\hline
\noalign{\smallskip}\hline\noalign{\smallskip}
\makecell[c]{$0$ \\ (Supple- \\ ment)} & \makecell[c]{$\footnotemark[1][u_{xx}, u_{xxxx}, uu_x,$ \\ $uu_{xxx}, uu_{xxxxx},$ \\ $u_{x}u_{xxxx}, u_{xx}u_{xxx},$ \\ $u_{xx}u_{xxxxx}, u_{xxx}u_{xxxx}]$} & \makecell[c]{$[u_{xx}, u_{xxxx}, uu_x]$} & \makecell[c]{$[uu_x]$}\\
 & (\textbf{\textcolor{blue}{-146,610.75}}) & (\textcolor{blue}{-135,102.20}) \checkmark & (\textcolor{blue}{-67,942.84})\\
\noalign{\smallskip}\hline
\end{tabular}
\addtolength{\tabcolsep}{3pt}
\end{center}
\footnotesize \footnotemark[1]{To avoid the minor details of cluttered discoveries, STRidge gets recursively reiterated with small magnitude coefficient removal until \(\forall j, \abs{\hat{\xi}_{j}}>10^{-1}\).}
\caption{\textbf{KS regularization hyperparameter selection}: STRidge's \(\mu\) is set to \((2(10^{2}), 5(10^{3}), 5(10^{3}))\), and \(d_{tol}\) equals \((1, 1, 50)\) for the three noise conditions. For the noisy \((x, t)\)\&\(u\) case, each polynomial candidate is normalized by its \(L_{1}\)-norm to get the better three-term PDE in terms of the BIC score.}
\vspace*{-\baselineskip}
\vspace*{-\baselineskip}
\label{tab:lambdas_ks}
\end{table}

We list the possible PDEs provided by STRidge for the various set of regularization hyperparameters in Table \ref{tab:lambdas_ks}. The metadata was specified as the 21,000 samples (\(N_{\mathcal{M}}\)) within the \([0, 100]\times[0,8]\) boundary generated by a Latin Hypercube Strategy \cite{stein1987large}. It alludes to us that the \(\lambda_{STR} = 10^{-5}\) founded PDEs cannot correspond to any of the \(\lambda_1 > 0\) feature importance because of the inclusion of \(u_{xxx}\), which may be inessential. Conversely, if we were to solely contemplate on the resulted PDEs associated with \(\lambda_1 = 0\), we would suspect that some terms are missing from \([u_{xx}, u_{xxxx}, uu_{x}]\) as the big PDE model comprising \([uu_{xxx}, uu_{xxxxx}, \dots, u_{xxx}u_{xxxx}]\) whose coefficient magnitudes were all comparable in size, e.g., of order \(> 10^{-1}\), demonstrated the lowest BIC score. The dilemma signifies that the unaided BIC, whose value varies dominantly by the changing log-likelihood term, cannot righteously balance the model complexity and accuracy, partly because no parsimonious governing PDE is involved behind the criterion assumption. In fact, the well-matched BIC is achievable by the simpler model built on the three correct candidates in \(\lambda_1 = 2(10^{-3})\). We mark the correct PDE expression $u_t = -0.989019u_{xx}-0.962360u_{xxxx}-0.966931uu_x$ found by \(\lambda_1=0\) as inferior to the selected model (\pmb{\checkmark}) in terms of discovery precision. \(\lambda_{STR} = 10^{-1}\) offers us the sparse PDEs, still, their BIC scores are much higher along with the clear BIC worthy enhancements observed when comparing against \(\lambda_{STR} = 10^{-3}\), thus designated as the condition giving the underfitting models. We take the PDE with the lowest BIC $u_t = -0.989305u_{xx}-0.970189u_{xxxx}-0.978123uu_x$ as our starting PDE (\pmb{\checkmark}), after assessing the agreed models for each \(\lambda_1 > 0\) row. For the noisy cases, the initial discovered KS \%CE are listed in Table \ref{tab:noise} (see the nPIML: IPI row). On the subsequent learning \circled{\textbf{(3)}} of Fig. \ref{fig:overview_framework}, the first (repolluted, if noisy) 21,504 data points were employed to finetune dPINNs.

In KS example, it is helpful to beware that including the higher-order derivatives in the basis candidates indicates enlarging the library size, which may have an ill effect on the discovery results. For example in the noisy \((x, t)\)\&\(u\) case, if we include \(u_{xxxxxx}\) and generate up to the 20-degree polynomials, the Pareto-optimal PDE with three terms, produced by \(\lambda_{STR}\)-varied STRidge, is wrong: $u_t = 0.524544u_{xx}-1.120505uu_x-0.626087uu_{xxx}$ (BIC = -93,841.66) instead of the previously found $u_t = -0.845746u_{xx}-0.818840u_{xxxx}-0.913990uu_x$ (BIC = -104,867.55). Nonetheless, this specific issue can be solved by searching over all possible PDEs with three terms to find the best PDE that shows the minimal BIC.

\begin{figure}
\centering
\includegraphics[width=0.5\textwidth]{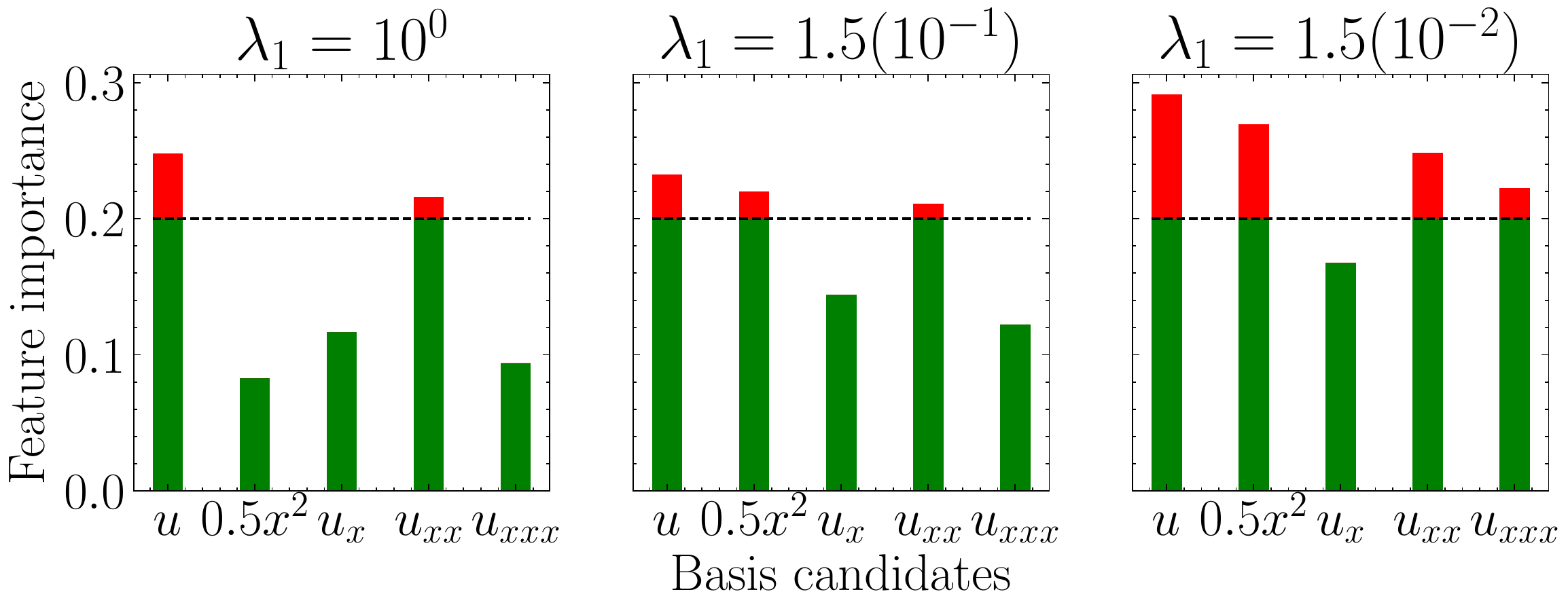}
\caption{QHO: Learned feature importance with varied \(\lambda_{1}\)}
\label{fig:feature_importance_qho}
\end{figure}

\begin{table}
\begin{center}
\begin{tabular}{cccc}
\hline\noalign{\smallskip}
$\lambda_{1}$/$\footnotemark[1]\lambda_{STR}$ & \footnotemark[2]$10^{-5}$ & $10^{-3}$ & $10^{-1}$\\
\noalign{\smallskip}\hline\noalign{\smallskip}
$10^{0}$ & \makecell[l]{$[u_x, u_{xx}, uu_x,$ \\ $\quad uu_{xx}, 0.5x^{2}u]$} & \makecell[c]{$[uu_{xx}, 0.5x^{2}u]$} & \makecell[c]{$[u]$}\\
 & (\textbf{-356,020.01}) & (-355,726.94) & (\textcolor{blue}{144,126.16})\\
\noalign{\smallskip}\hline\noalign{\smallskip}
$1.5(10^{-1})$ & \makecell[l]{$[u_x, u_{xx},$ \\ $\quad uu_{x}, 0.5x^{2}u]$} & \makecell[c]{$[u_{xx}, 0.5x^{2}u]$} & \makecell[c]{$[u]$}\\
 & (\textbf{-325,329.72}) & (\textcolor{blue}{-325,110.45}) \checkmark & (\textcolor{blue}{144,190.29})\\
\noalign{\smallskip}\hline\noalign{\smallskip}
$1.5(10^{-2})$ & \makecell[l]{$[u_x, u_{xx},$ \\ $\quad uu_{x}, 0.5x^{2}u]$} & \makecell[c]{$[u_{xx}, 0.5x^{2}u]$} & \makecell[c]{$[u]$}\\
 & (\textbf{-325,526.78}) & (\textcolor{blue}{-325,307.53}) \pmb{\checkmark} & (\textcolor{blue}{144,195.08})\\
\noalign{\smallskip}\hline
\end{tabular}
\end{center}
\footnotesize \footnotemark[1]{We enumerate \(\lambda_{STR}\) from \((10^{-3}, 10^{-2}, 10^{-1})\) for the noisy \((x, t)\)\&\(u\) case.}
\footnotesize \footnotemark[2]{STRidge is refitted once to show only the term that \(\abs{\hat{\xi}_{j}}>1.4(10^{-2})\).}
\caption{\textbf{QHO regularization hyperparameter selection}: STRidge's \(\lambda_{0}\) is set to \(10^{2}\lambda_{STR}\varepsilon\), and \(d_{tol}\) equals \(10\) for the three noise conditions.}
\vspace*{-\baselineskip}
\vspace*{-\baselineskip}
\label{tab:lambdas_qho}
\end{table}

\subsubsection{Initial Discovered QHO PDE} \label{discovering_qho}
As per the specific treatments for QHO mentioned in \ref{exp_settings}, \textbf{Algorithm \ref{alg:1}} turns applicable for the complex-valued PDEs. The preselector was trained with varied \(\lambda_1\). Each basis candidate importance at the different levels is shown in Fig. \ref{fig:feature_importance_qho}. All three correct candidates can surpass the threshold when \(\lambda_1=1.5(10^{-1})\) or \(1.5(10^{-2})\) whereas \(\lambda_1=10^{0}\) compels the too strong regularization.

For QHO, the metadata for STRidge was the linearly discretized points from the full-field spatio-temporal domain, i.e, \(N_{\mathcal{M}}=82,432\) and \(\forall i, t^{\mathcal{M}}_{i} \leqslant 4\). The cross results for the regularization hyperparameter selection are listed in Table \ref{tab:lambdas_qho}. If the \(\lambda_{STR}\) intensity is loosen from \(10^{-1}\) to \(10^{-3}\) the considerable shoots in the BIC improvement are apparently gained. However, regularizing too mildly, e.g., \(\lambda_{STR}=10^{-5}\), does not provide any left necessary candidates, exhibiting the small BIC reductions with the more unsound terms that are unstable across varying \(\lambda_1\). Ultimately, $u_t = (-0.000463+0.498906i)u_{xx}+(-0.002272-0.999284i)0.5x^{2}u$ (\pmb{\checkmark}) is accepted for the denoising and finetuning stage owing to its minimal BIC score among the agreed PDEs. In the cases where noise exists, the \%CE of the initial discovered complex-valued PDEs are shown in Table \ref{tab:noise} (see the nPIML: IPI row).

\begin{figure}
\centering
\includegraphics[width=0.5\textwidth]{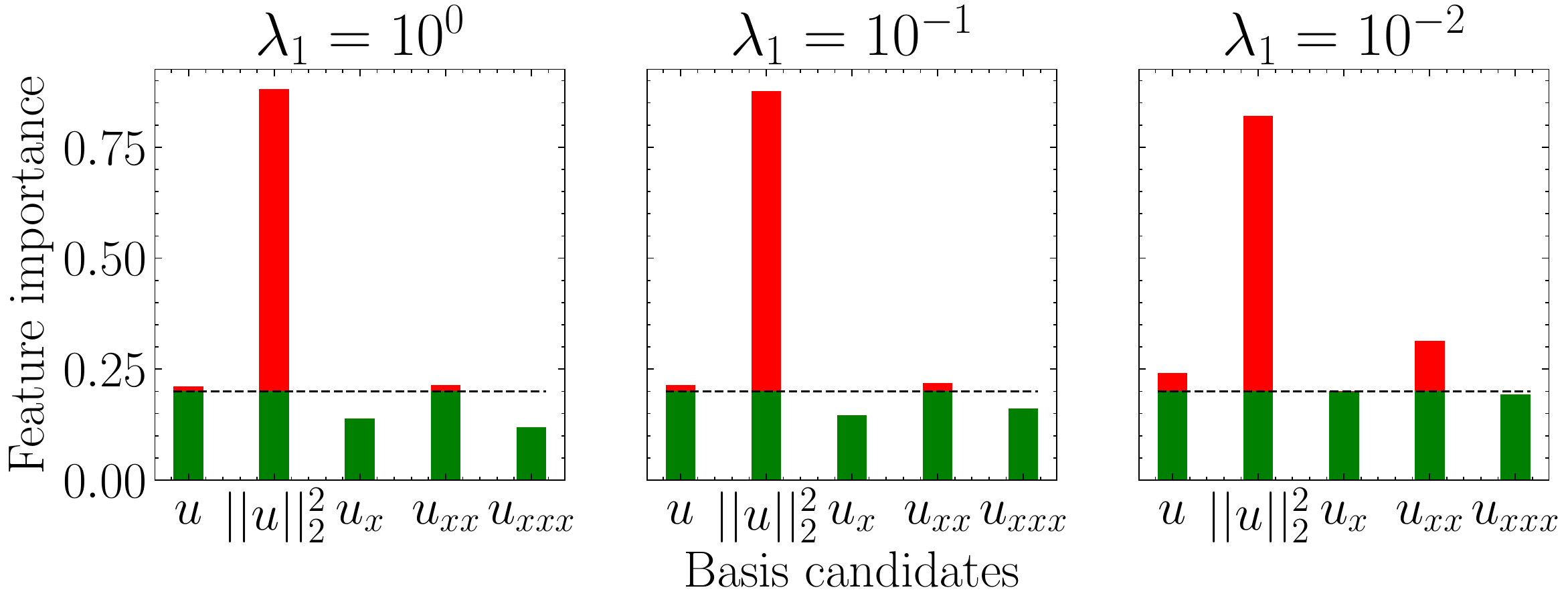}
\caption{NLS: Learned feature importance with varied \(\lambda_{1}\)}
\label{fig:feature_importance_nls}
\end{figure}

\begin{table}
\begin{center}
\setlength{\tabcolsep}{3pt}
\begin{tabular}{cccc}
\hline\noalign{\smallskip}
$\lambda_{1}$/$\lambda_{STR}$ & \footnotemark[1]$10^{-7}$ & $10^{-5}$ & $10^{-2}$\\
\noalign{\smallskip}\hline\noalign{\smallskip}
$10^{0}$ & \makecell[c]{$[u, \norm{u}^{2}_{2}, u_x, u_{xx},$ \\ $\quad u^2, u\norm{u}^{2}_{2}, uu_x, \norm{u}^{2}_{2}u_x]$} & \makecell[c]{$[u_{xx}, u\norm{u}^{2}_{2}]$} & \makecell[c]{$[u\norm{u}^{2}_{2}]$}\\
 & (\textbf{-108,572.98}) & (\textcolor{blue}{-107,799.25}) \checkmark & (\textcolor{blue}{166,786.21})\\
\noalign{\smallskip}\hline\noalign{\smallskip}
$10^{-1}$ & \makecell[c]{$[u_x, u_{xx}, u^2, u\norm{u}^{2}_{2},$ \\ $\quad uu_x, uu_{xx}, u^{2}_{x}]$} & \makecell[c]{$[u_{xx}, u\norm{u}^{2}_{2}]$} & \makecell[c]{$[u\norm{u}^{2}_{2}]$}\\
 & (\textbf{-108,582.11}) & (\textcolor{blue}{-107,847.00}) \pmb{\checkmark} & (\textcolor{blue}{166,790.61})\\
\noalign{\smallskip}\hline\noalign{\smallskip}
$10^{-2}$ & \makecell[c]{$\footnotemark[2][u_x, u_{xx}, u^2, u\norm{u}^{2}_{2},$ \\ $\quad uu_x, uu_{xx}, u^{2}_{x}]$} & \makecell[c]{$[u_{xx}, u\norm{u}^{2}_{2}]$} & \makecell[c]{$[u\norm{u}^{2}_{2}]$}\\
 & (\textbf{-108,508.42}) & (\textcolor{blue}{-107,766.91}) \checkmark & (\textcolor{blue}{166,790.63})\\
\noalign{\smallskip}\hline
\end{tabular}
\end{center}
\footnotesize \footnotemark[1]{STRidge is refitted once to show only the term that \(\abs{\hat{\xi}_{j}}>1.4(10^{-3})\).}
\footnotesize \footnotemark[2]{\((0.000294-0.000829i)u_{xxx}\) that partly causes the disagreement is withdrawn from the list, since \(\abs{0.000294-0.000829i} \leqslant 1.4(10^{-3})\).}
\caption{\textbf{NLS Regularization hyperparameter selection}: STRidge's \(\lambda_{0}\) is set to \(10^{5}\lambda_{STR}\varepsilon\), and \(d_{tol}\) equals \(100\) for the three noise conditions.}
\vspace*{-\baselineskip}
\vspace*{-\baselineskip}
\label{tab:lambdas_nls}
\end{table}

\subsubsection{Initial Discovered NLS PDE} \label{discovering_nls}
The feature importance measures are displayed in Fig. \ref{fig:feature_importance_nls}. The correct candidates are safely secured, passing the threshold and becoming effective for all the choices of \(\lambda_{1}=10^{0}\), \(10^{-1}\) or \(10^{-2}\). Despite that, \(\lambda_1=10^{-2}\) is relatively low such that the inclusion of \(u_x\) might have complicated the hyperparameter selection procedure.

We limited the whole domain arbitrarily at \(t<1.25\) for bounding the interested region upon which the metadata was linearly discretized, i.e., in total \(N_{\mathcal{M}}=40,960\). Still, we positively ensured that the essential dynamics were covered. The found PDEs are assimilated in Table \ref{tab:lambdas_nls}, indexing diverse set of \((\lambda_1, \lambda_{STR})\) for the regularization hyperparameter selection. The admittance of \(u_{xx}\), resulted by decreasing \(\lambda_{STR}\) from \(10^{-2}\) to \(10^{-5}\), apparently upgrades the BIC scores. Further dropping \(\lambda_{STR}\) down to \(10^{-7}\) can push the BIC scores down slightly with the increased terms that, however, end up disagreeing with the preselectors. Like QHO example in \ref{discovering_qho}, the agreed sparse PDE, exhibiting the minimal BIC, gets accepted to be denoised and finetuned. For NLS, the initial discovered PDE reads $u_t = (-0.000863+0.499928i)u_{xx}+(-0.000973+0.999259i)u\norm{u}^{2}_{2}$ (\pmb{\checkmark}). In the noisy experiments, the \%CE of the initial discovered complex-valued PDEs are provided in Table \ref{tab:noise} (see the nPIML: IPI row).

\begin{table*}[t]
\begin{center}
\resizebox{\textwidth}{!}{\begin{tabular}{llllll}
\hline\noalign{\smallskip}
Dataset & Method & \# Train samples (\(N_f\)) & Noiseless & \(u\) + Noise\textsubscript{\(u\)} & \(u\) + Noise\textsubscript{\(u\)} \& \((x, t)\) + Noise\textsubscript{\((x, t)\)}\\
\noalign{\smallskip}\hline\noalign{\smallskip}
   & PDE-FIND (STRidge) \cite{rudy2017data} & 256$\times$100\footnotemark[1] & 19.2070$\pm$19.0686 & Failed ($-0.0698uu_{x}$) & Not applicable\footnotemark[2]\\
   & DLrSR \cite{li2020robust} & 256$\times$100 & 19.2070$\pm$19.0686 & Failed ($-0.0698uu_{x}$) & Not applicable\\
Burgers & PINN\footnotemark[3] \cite{raissi2019physics} & 3,000 & 0.3256$\pm$0.1921 & 0.9212$\pm$0.8589 & 4.0893$\pm$2.9622\\
   & nPIML: IPI\footnotemark[4] & 3,000 & 2.5730$\pm$1.1904 & 7.0093$\pm$2.6069 & 55.2051$\pm$15.8919\\
   & nPIML w/o Denoise\footnotemark[5] & 3,000 & 0.1264$\pm$0.0605 & 0.4271$\pm$0.2451 & 2.9920$\pm$2.2222\\
   & nPIML & 3,000 & \textbf{0.0557$\pm$0.0170} & \textbf{0.3360$\pm$0.1251} & \textbf{0.8546$\pm$0.4806}\\
\noalign{\smallskip}\hline
   & PDE-FIND (STRidge) & 128$\times$501 & 0.5194$\pm$0.1733 & Failed ($-5.4128uu_{x}$) & Not applicable\\
   & DLrSR & 128$\times$501 & 0.5194$\pm$0.1733 & Failed ($-5.3521uu_{x}$) & Not applicable\\
KdV & nPIML: IPI & 2,000 & 0.8710$\pm$0.2224 & 2.9887$\pm$1.1612\footnotemark[6] & 3.7460$\pm$1.4158\footnotemark[6]\\
   & nPIML w/o Denoise & 2,000 & 0.6413$\pm$0.3904 & 1.2547$\pm$0.8369 & 2.9378$\pm$1.6140\\
   & nPIML & 2,000 & \textbf{0.0890$\pm$0.0568} & \textbf{0.2845$\pm$0.2463} & \textbf{0.4344$\pm$0.2696}\\
\noalign{\smallskip}\hline
   & PDE-FIND (STRidge) & 1024$\times$251 & 0.7557$\pm$0.5967 & 52.2843$\pm$1.4005 & Not applicable\\
   & DLrSR & 1024$\times$251 & 0.7571$\pm$0.5966 & Failed\footnotemark[7] & Not applicable\\
KS & nPIML: IPI & 80,000 & 2.0794$\pm$0.7842 & 10.7558$\pm$3.3449 & 14.0475$\pm$4.0048\\
   & nPIML w/o Denoise & 80,000 & 1.7417$\pm$1.1171 & 8.8925$\pm$5.2704 & 9.2365$\pm$6.5974\\
   & nPIML & 80,000 & \textbf{0.4775$\pm$0.2751} & \textbf{2.9320$\pm$1.4401} & \textbf{3.6493$\pm$3.9688}\\
\noalign{\smallskip}\hline
   & PDE-FIND (STRidge) & 512$\times$161 & 0.2458$\pm$0.0101 & 9.3850$\pm$6.7242 & Not applicable\\
   & DLrSR & 512$\times$161 & 0.2850$\pm$0.0090 & 9.3711$\pm$6.7143 & Not applicable\\
QHO & nPIML: IPI & 30,000 & 0.2379$\pm$0.0003 & 0.3163$\pm$0.0705 & 0.4197$\pm$0.0121\\
   & nPIML w/o Denoise & 30,000 & 0.0377$\pm$0.0211 & 0.2380$\pm$0.1463 & 0.3278$\pm$0.1694\\
   & nPIML & 30,000 & \textbf{0.0278$\pm$0.0193} & \textbf{0.1235$\pm$0.0580} & \textbf{0.2669$\pm$0.1639}\\
\noalign{\smallskip}\hline
   & PDE-FIND (STRidge) & 256$\times$201 & 0.3469$\pm$0.2888 & 2.8485$\pm$2.6764 & Not applicable\\
   & DLrSR & 256$\times$201 & 0.3294$\pm$0.2801 & 2.8542$\pm$2.6778 & Not applicable\\
NLS & nPIML: IPI & 2,500 & 0.1478$\pm$0.0255 & 0.5686$\pm$0.2517 & 2.3726$\pm$1.5939\\
   & nPIML w/o Denoise & 2,500 & 0.0491$\pm$0.0060 & 0.0953$\pm$0.0114 & 0.2205$\pm$0.0877\\
   & nPIML & 2,500 & \textbf{0.0421$\pm$0.0172} & \textbf{0.0571$\pm$0.01327} &
   \textbf{0.1652$\pm$0.0532}\\
\noalign{\smallskip}\hline
\end{tabular}}
\end{center}
\footnotesize \footnotemark[1]{All the discretized points are shown in the mesh representation: \# in \(x \times t\).}
\footnotesize \footnotemark[2]{Because a mesh is required for taking polynomial derivatives used in PDE-FIND where STRidge was firstly introduced.}
\footnotesize \footnotemark[3]{\(\hat{\xi}\) is initialized at \(\mqty[exp(-7.0), 1.0]^{\intercal}\) before training PINN.}
\footnotesize \footnotemark[4]{The results until \circled{\textbf{(2)}} of Fig. \ref{fig:overview_framework}, Initial PDE Identification.}
\footnotesize \footnotemark[5]{The results from \circled{\textbf{(3)}} of Fig. \ref{fig:overview_framework}, dPINNs, but without the denoising DFT module and projection networks.}
\footnotesize \footnotemark[6]{\(\lambda_1\) is assigned to \(2(10^{-5})\) instead of \(2(10^{-6})\).}
\footnotesize \footnotemark[7]{DLrSR with the original and unvarying \(\lambda_1\) discovers the following mismatched PDE: \(u_t=-0.60uu_{x}-0.39u_{xx}-0.10uu_{xxx}-0.49u_{xxxx}\).}
\caption{\textbf{Summary of the robust discovery results by nPIML}: The noise is 1\% of standard deviation. Generally, the adopted \(\lambda_1\)s for the noisy experiments are identical to the noiseless condition unless noted otherwise. The best error is \textbf{bolded}.}
\vspace*{-\baselineskip}
\label{tab:noise}
\end{table*}

\subsection{Finetuning PDE Coefficients by dPINNs}
Based on the results in Table \ref{tab:noise}, nPIML establishes superior results over nPIML without the denoising DFT and projection networks for the noisy cases, especially when both \((x, t)\) and \(u\) are contaminated. For the clean dataset, the denoising mechanism seems to not over perturb backwardly through converging \(\beta_{(x, t)}, \beta_{u} \rightarrow 0\), maintaining the effectiveness of the dPINNs' learning by \textbf{Algorithm \ref{alg:2}}, on par to the nPIML without the denoising that exactly matches the noiseless hypothesis. Indeed, nPIML can outperform nPIML without the denoisers since the shifting to the more propitious finite set, e.g., \(\set{(x^{*}_{i}, t^{*}_{i}, u^{*}_{i})^{N_f}_{i=1}}\), is still technically probable. In Burgers' example, nPIML surpasses vanilla PINN for all experimental cases regardless of the denoising modules, implying the superiority and benefits of the precomputed initialization followed by finetuning \(\hat{\theta}\) and \(\hat{\xi}\). Moreover, if the genuine PDE is known beforehand, training PINN from scratch eventually leads to the better close-formed discovery than PDE-FIND (STRidge). The accuracy enhancement points out the usefulness of automatic differentiation and physics-informed learning.

\subsection{Robustness against Scarce Data} \label{robust_against_scarce_data}
Table \ref{tab:tolerance} reveals the tolerance against the decreasing number of training samples in Burgers' example. The precise discovered PDEs are obtainable by finetuning the coefficients even though only the 500 training data points are available. However, it is challenging to recover Burgers' PDE if the noise is added or dPINNs are trained with just the 100 training samples, implied by the faulty discoveries by \textbf{Algorithm \ref{alg:1}}. Fortunately, the results show that the pragmatic denoising affine transformation by the projection networks is feasible even under the noisy and moderately limited number of labeled samples, e.g., 1,000. It is worth pointing out that data bias towards diverse training sets leads to diversity in (initial) discovery results when learning from a few samples. In addition, the involving parameter and model initializations affect PINN approximated outputs as discussed in \cite{wong2021learning}, and undoubtedly the PDEs that are derived from those outputs.

\begin{table}
\begin{center}
\resizebox{0.48925\textwidth}{!}{
\begin{tabular}{llll|l}
\hline\noalign{\smallskip}
\# Train & Noise & nPIML: IPI \%CE & \multicolumn{2}{c}{Finetuned PDE \%CE}\\
samples & & & w/o Denoise & w/ Denoise\\
\noalign{\smallskip}\hline\noalign{\smallskip}
3000\footnotemark[1] & N\footnotemark[2] & 2.5730$\pm$1.1904 & 0.1264$\pm$0.0605 & \textbf{0.0557$\pm$0.0170}\\
 & Y\footnotemark[3] & 55.2051$\pm$15.8919 & 2.9920$\pm$2.2222 & \textbf{0.8546$\pm$0.4806}\\
\noalign{\smallskip}\hline\noalign{\smallskip}
1000 & N & 3.8530$\pm$1.6829 & 1.0953$\pm$1.0526 & \textbf{0.8105$\pm$0.7565}\\
 & Y & 26.5837$\pm$2.6611 & 8.3633$\pm$8.2076 & \textbf{1.6114$\pm$1.1907}\\
\noalign{\smallskip}\hline\noalign{\smallskip} 
500 & N & 6.2883$\pm$2.6029 & 1.9302$\pm$1.7908 & \textbf{1.4888$\pm$1.2651}\\
 & Y & \makecell[c]{Failed\footnotemark[4]} & \makecell[c]{Not applicable} & \makecell[c]{Not applicable}\\
\noalign{\smallskip}\hline\noalign{\smallskip} 
100 & N & \makecell[c]{Failed\footnotemark[5]} & \makecell[c]{Not applicable} & \makecell[c]{Not applicable}\\
 & Y & \makecell[c]{Failed\footnotemark[6]} & \makecell[c]{Not applicable} & \makecell[c]{Not applicable}\\
\noalign{\smallskip}\hline
\end{tabular}
}
\end{center}
\footnotesize \footnotemark[1]{Taken from Table \ref{tab:noise}.}
\footnotesize \footnotemark[2]{Noiseless.}
\footnotesize \footnotemark[3]{\(u\) + Noise\textsubscript{\(u\)} \& \((x, t)\) + Noise\textsubscript{\((x, t)\)}.\\}
\footnotesize \footnotemark[4]{\(u_t=-0.703435uu_{x}-0.000041u_{x}u_{xx}\).}
\footnotesize \footnotemark[5]{\(u_t=-0.509967uu_{x}\).\\}
\footnotesize \footnotemark[6]{\(u_t=-0.621560uu_{x}\).}
\caption{\textbf{Discovered Burgers' PDE on the scarce data}}
\vspace*{-\baselineskip}
\vspace*{-\baselineskip}
\label{tab:tolerance}
\end{table}

\begin{figure}
\begin{subfigure}[b]{0.45\textwidth}
  \centering
  \includegraphics[width=\textwidth]{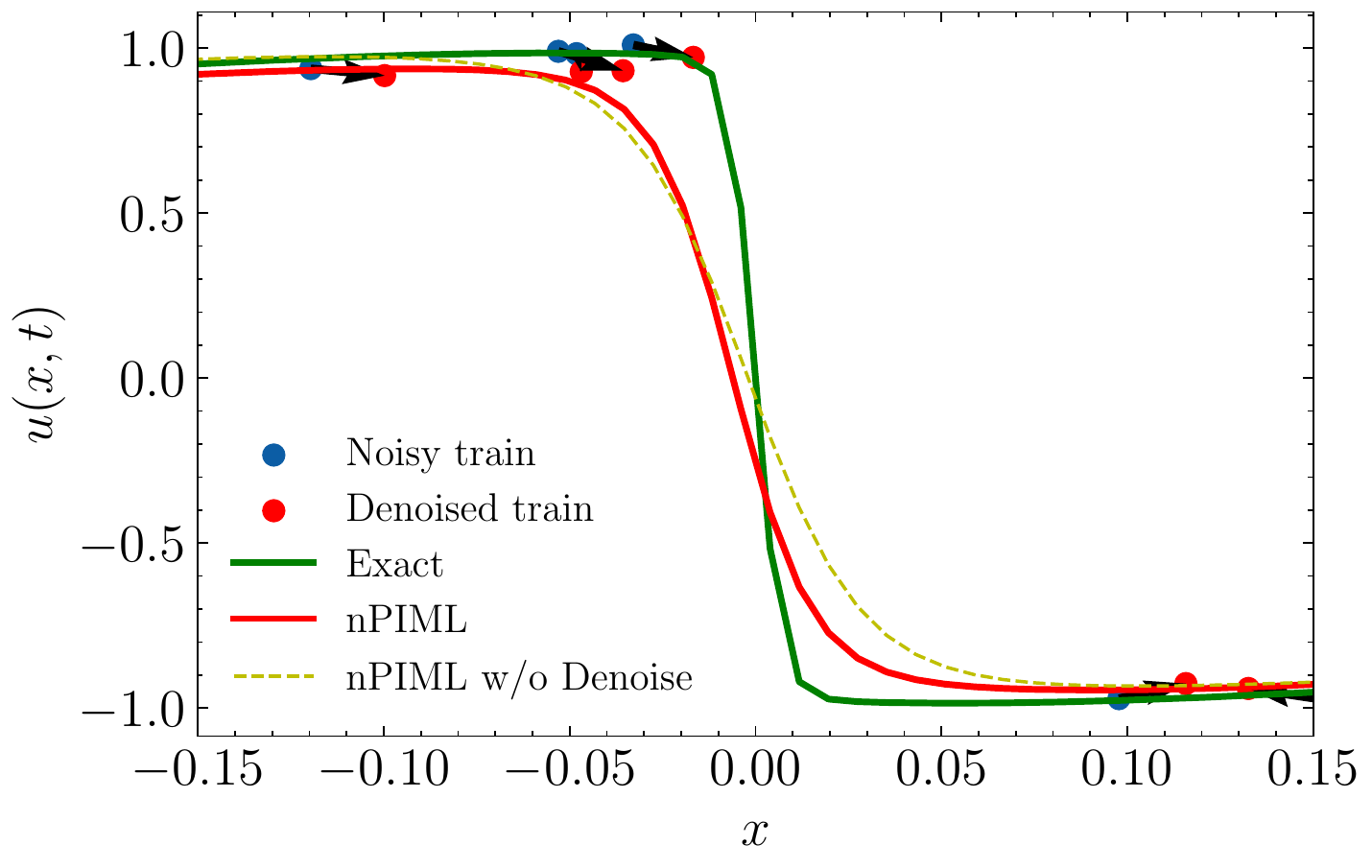}
  \caption{Snapshot around \(t=0.46\)}
  \label{fig:dft_vs_without_dft_t=046}
\end{subfigure}
\newline
\begin{subfigure}[b]{0.45\textwidth}
  \centering
  \includegraphics[width=\textwidth]{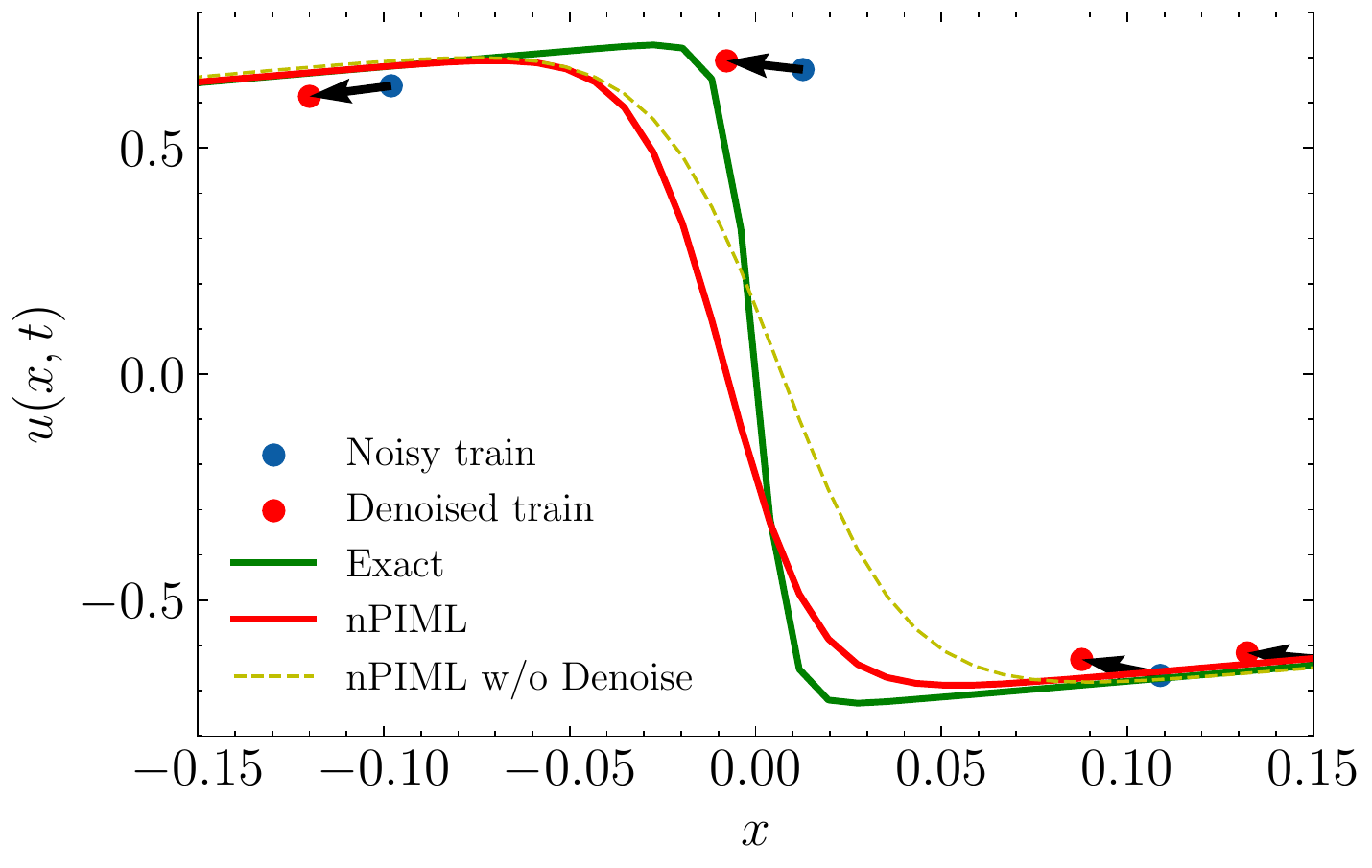}
  \caption{Snapshot around \(t=0.97\)}
  \label{fig:dft_vs_without_dft_zoom_t=097}
\end{subfigure}
\caption{Close visualization of how the preselector networks react to the high noise at \(x \in [-0.15, 0.15]\), around the abrupt transition caused by the shock waves.}
\vspace*{-\baselineskip}
\label{fig:dft_vs_without_dft_zoom}
\end{figure}

\subsection{Denoising Mechanism against High Noise}
\subsubsection{Denoising Visualization} \label{denoising_visualization}
We sought to apprehend how the projection networks respond to high noise visually by letting dPINNs expose the strongly contaminated dataset, where \(u\) and \((x, t)\) are polluted with \(noise(u, 5)\) and \(noise((x, t), 5)\). Specifically, we finetuned the dPINNs pretrained by \(\hat{\theta}\), taken from the 1\%Noise\(+(x, t)\)\&\(u\) case of Burgers' PDE. The initialized PDE was resolved by LS based on the intentionally uplifted 5\% noisy \((x, t)\), expressing the form as follows: $u_t = 0.000606u_{xx}-0.403049uu_{x}$. For such high noise, we find it is useful that \(\mathcal{P}_{\Omega_{(x, t)}}(x, t)\) and \(\mathcal{P}_{\Omega_{u}}(u)\) should not be only activated by the final Tanh but also unbiased standardized and then scaled down to be 0.01 times the values to denoise gradually from small to larger noise magnitude since denoising the considerable amount at the beginning of the dPINNs' learning can ultimately cause the divergence. \(\alpha\) and \((\beta^{\prime}_{(x, t)}, \beta^{\prime}_{u})\) are initialized at 0.1 and \((10^{-3}, 10^{-3})\). We display how the projection networks denoise closely around \(t=0.46, 0.97\) in Fig. \ref{fig:dft_vs_without_dft_zoom}. By the proximate examination near the dynamically changing region, where there are only a few supervised samples, the naive PDE estimation neglecting the noise effect is observed if the denoising components are ablated. The optimized PDE is $u_t = 0.012378u_{xx}-0.948156uu_{x}$ without the denoising. In comparison, the projection networks can shift the polluted samples towards the direction that drives the approximated solution by dPINNs to better captures the exact characteristics of Burgers' PDE when the denoising components are utilized. As seen in Fig. \ref{fig:dft_vs_without_dft_t=046} and \ref{fig:dft_vs_without_dft_zoom_t=097}, the noisy samples that are more to the right of the exact solution, get redirected to the left and vice versa. Positively impacted by the denoising, the optimized PDE carries the form of $u_t = 0.008550u_{xx}-0.972390uu_{x}$.

\subsubsection{Finetuning against High Noise} \label{finetuning_high_noise}
Since restoring a decent approximation of the hidden KS PDE from highly noisy data can be sensitive and challenging. We, therefore, set up more experiments, similar to \ref{denoising_visualization}, finetuning the dPINNs initialized with \(\hat{\theta}\) taken from the 1\%Noise\(+(x, t)\)\&\(u\) case, but single \(\hat{\xi}\) uniformly generated such that \(\forall i, \hat{\xi}_{i}\) \(\sim (-10^{-6})\boldsymbol{\mathcal{U}}(0, 1)\). The intensity of the noise that contaminates \((x, t)\)\&\(u\) is explicitly increased to 3\%, 5\% and 10\%. \(\alpha\) and \((\beta^{\prime}_{(x, t)}, \beta^{\prime}_{u})\) are initialized at \(0.1\) and \((10^{-2}, 10^{-2})\). \(\beta^{\prime}_{(x, t)}\) and \(\beta^{\prime}_{u}\) are clamped within \([-1.0, 1.0]\) during the finetuning process. The quantitative results in Table \ref{tab:high_noise}, even more, emphasize the superiority of asserting the denoising mechanism to minimize the discovery error numerically under much-corrupted dataset.

\begin{table}
\begin{center}
\begin{tabular}{cc|c}
\hline\noalign{\smallskip}
Noise level & \multicolumn{2}{c}{Finetuned PDE \%CE}\\
 & w/o Denoise & w/ Denoise\\
\noalign{\smallskip}\hline\noalign{\smallskip}
1\%\footnotemark[1] & 9.2365$\pm$6.5974 & \textbf{3.6493$\pm$3.9688}\\
3\% & 27.0814$\pm$20.0158 & \textbf{11.5509$\pm$9.7542}\\
5\% & 45.8996$\pm$31.3289 & \textbf{20.5851$\pm$24.1741}\\
10\% & 56.8900$\pm$40.8643 & \textbf{52.3366$\pm$38.7182}\\
\hline\noalign{\smallskip}
\footnotesize \footnotemark[1]{Taken from Table \ref{tab:noise}.}
\end{tabular}
\end{center}
\caption{\textbf{Numerical results of finetuning dPINNs on highly noisy KS data}}
\vspace*{-\baselineskip}
\vspace*{-\baselineskip}
\label{tab:high_noise}
\end{table}

\section{Conclusion} \label{conclusion}
We have presented the interpretable and noise-aware physics-informed machine learning framework for distilling the nonlinear PDE governing a physical system in an analytical expression. The proposed method mainly tackles the problems with the suboptimal derivatives, sensitivity of regularization hyperparameters, and polluted datasets. The weakly physics-informed solver network is the primary building block for derivative computation. Multi-perspective assessment of the diverse sets of regularization hyperparameters is feasible through the physics-learning preselector network and the sparse regression. Finally, denoising physics-informed neural networks are introduced to finetune the objective PDE coefficients to the optimality on the affine transformed noise-reduced dataset given by the projection networks. The numerical results show that the proposed method is robust to the scarcity of labeled samples and noise on five classic canonical PDEs, outperforming the state-of-art regression-based discovery methods.

Nonetheless, the proposed framework exhibits some limitations. For instance, there is no explicit denoising mechanism at the early derivative preparation and sparse regression stages; thus, particular noise of an unknown distribution may fake those initial processes and let the entire framework fail. The predicament that underlying physics remains mysterious initially causes the projection networks to be inoperable as the affine transformation can yield the unwanted \(\Tilde{u} \approx \vv{0}\), and solely assigning an appropriate threshold for denoising DFT is not either trivial or readily beneficial. Towards future improvements, researchers may conduct extensive studies on grounded topics such as the effect of parameter initialization on the discovery stability or a border class of inferable PDEs that is not restricted by the linear assumption.

{\small
\bibliographystyle{ieeetr}
\bibliography{ref}
}

\end{document}